\documentstyle[amssymb,newlfont]{amsart}

\newtheorem{thm}{Theorem}[section]
\newtheorem{lem}[thm]{Lemma}
\newtheorem{prop}[thm]{Proposition}
\newtheorem{cor}[thm]{Corollary}

\theoremstyle{definition}
\newtheorem{defn}{Definition}

\theoremstyle{remark}
\newtheorem{notation}{Notation}

\newtheorem{rem}{Remark}

\numberwithin{equation}{section}

\newcommand{\C}{\Bbb{C}}

\newcommand{\Q}{{\Bbb{Q}}}
\newcommand{\Z}{\Bbb{Z}}
\newcommand{\N}{\Bbb{N}}

\newcommand{\K}{\Bbb{K}}

\def\Ker{\operatorname {Ker}}

\def\End{\operatorname {End}}

\def\codim{\operatorname {codim}}
\def\Prim{\operatorname {Prim}}
\def\Spec{\operatorname {Spec}}
\def\Symp{\operatorname {Symp}}
\def\Lie{\operatorname {Lie}}
\def\ad{\operatorname {ad}}
\def\rk{\operatorname {rk}}
\def\Ad{\operatorname {Ad}}

\def\Aut{\operatorname {Aut}}
\def\Soc{\operatorname {Soc}}
\def\Ann{\operatorname {Ann}}
\def\Stab{\operatorname {Stab}}

\newcommand{\smatrix}[4]{\bigl[ \begin{smallmatrix} {#1} & {#2} \\ {#3}
& {#4} \end{smallmatrix} \bigr]}
\newcommand{\isomto}{\, \overset{\sim}{\rightarrow} \, }

\newcommand{\pol}[1]{{#1}^\circ}
\newcommand{\orth}[1]{{#1}^\perp}
\newcommand{\pair}[2]{\langle {#1},{#2} \rangle}
\newcommand{\rkZ}{\rk_{\Bbb Z} \, }
\newcommand{\hpair}[2]{\langle {#1} \mid {#2} \rangle}

\newcommand{\m}{\frak{m}}

\newcommand{\g}{\frak{g}}

\newcommand{\npm}{\frak{n}^{\pm}}

\newcommand{\tgoth}{\frak{t}}
\newcommand{\h}{\frak{h}}
\newcommand{\k}{\frak{k}}
\newcommand{\bm}{\frak{b}^-}
\newcommand{\bp}{\frak{b}^+}
\newcommand{\bpm}{\frak{b}^{\pm}}

\newcommand{\up}{\frak{u}^+}
\newcommand{\upm}{\frak{u}^{\pm}}

\newcommand{\agoth}{\frak{a}}
\newcommand{\atilde}{\tilde{\frak{a}}}
\newcommand{\hQ}{\frak{h}_{\Bbb Q}}
\newcommand{\hQdual}{\frak{h}^*_{\Bbb Q}}
\newcommand{\dgoth}{\frak{d}}

\newcommand{\bP}{\bold{P}}

\newcommand{\bQ}{\bold{Q}}
\newcommand{\bB}{\bold{B}}
\newcommand{\bK}{\bold{K}}
\newcommand{\bR}{\bold{R}}
\newcommand{\bRp}{\bold{R}_+}

\newcommand{\Char}[1]{\bold{X}({#1})}

\newcommand{\bL}{\bold{L}}

\newcommand{\Gbar}{\bar{G}}

\newcommand{\vart}{\vartheta}

\newcommand{\Phipm}{\Phi_{\pm}}
\newcommand{\Phip}{\Phi_+}
\newcommand{\Phim}{\Phi_-}

\newcommand{\wplus}{w_+}
\newcommand{\wminus}{w_-}
\newcommand{\wdot}{\dot{w}}
\newcommand{\Cw}{\cal{C}_w}
\newcommand{\Cwd}{\cal{C}_{\dot{w}}}
\newcommand{\Bw}{\cal{B}_w}
\newcommand{\Bwd}{\cal{B}_{\dot{w}}}
\newcommand{\Aw}{\cal{A}_w}
\newcommand{\Awd}{\cal{A}_{\dot{w}}}

\newcommand{\Uq}{U_q(\g)}
\newcommand{\Dq}{D_q(\g)}
\newcommand{\Dqp}{D_{q,p^{-1}}(\g)}
\newcommand{\uqbm}{U_q(\frak{b}^-)}
\newcommand{\uqpbm}{U_{q,p^{-1}}(\frak{b}^-)}
\newcommand{\uqbp}{U_q(\frak{b}^+)}
\newcommand{\uqpbp}{U_{q,p^{-1}}(\frak{b}^+)}
\newcommand{\uqnm}{U_q(\frak{n}^-)}
\newcommand{\uqnp}{U_q(\frak{n}^+)}
\newcommand{\Cq}{{\Bbb{C}}_q[G]}
\newcommand{\Cqp}{{\Bbb{C}}_{q,p}[G]}

\newcommand{\qexp}[1]{q^{{#1}}}
\newcommand{\cfv}[2]{c_{#1,#2}}
\newcommand{\ct}{\tilde{c}}
\newcommand{\ctw}[1]{\tilde{c}_{w{#1}}}
\newcommand{\cw}[1]{c_{w{#1}}}
 
\newcommand{\cqp}{{\Bbb{C}}_{q,p}[G]}

\newcommand{\ptil}{\tilde{p}}

\newcommand{\Aop}{A^{\text{op}}}
\newcommand{\calCq}{\cal{C}_q}
\newcommand{\calCqp}{\cal{C}_{q,p}}
\newcommand{\Dqpp}{D_{q,p}(\g)}
\newcommand{\cte}{\frac{4i\pi}{\hbar}}

\begin{document}


\title{Algebraic structure of multi-parameter quantum groups}
\author{Timothy J. Hodges} 
\thanks{The first author was partially supported by grants from the National 
Security Agency and the C. P. Taft Memorial Fund.}
\address{University of Cincinnati, Cincinnati, OH 45221-0025, U.S.A.}
\email{timothy.hodges@@uc.edu}
\author{Thierry Levasseur}
\address{Universit\'e de Poitiers, 86022 Poitiers,  France}
\email{thierry.levasseur@@matpts.univ-poitiers.fr}
\author{Margarita Toro}
\thanks{The third author was partially supported by a grant from Colciencias.}
\address{Universidad Nacional de Colombia, Apartado A\'{e}reo, Medellin, Colombia}
\email{mmtoro@@perseus.unalmed.edu.co}
\maketitle

\begin{center}	{\bf Introduction} 			\end{center}

\smallskip

	Let $G$ be a connected semi-simple complex Lie group. We
define and study the multi-parameter quantum group $\Cqp$ in the case where $q$ is a complex parameter that is not a root of unity. Using a
method of twisting bigraded Hopf algebras by a cocycle, \cite{AST}, we
develop a unified approach to the construction of $\Cqp$ and of the
multi-parameter Drinfeld double $D_{q,p}$. Using a general method of
deforming bigraded pairs of Hopf algebras, we construct a Hopf pairing
between these algebras from which we deduce a Peter-Weyl-type theorem
for $\cqp$. We then describe the prime and primitive spectra of
$\cqp$, generalizing a result of Joseph. In the one-parameter case
this description was conjectured, and established in the $SL(n)$-case,
by the first and second authors in
\cite{HL1,HL2}. It was proved in the general case by Joseph in
\cite{J1,J2}. In particular the
orbits in $\Prim \cqp$ under the natural action of the maximal torus
$H$ are indexed, as in the one-parameter case by the elements of the
double Weyl group $W \times W$.  
Unlike the one-parameter case there is not in general a
bijection between $\Symp G$ and $ \Prim \cqp$. However in the case
when the symplectic leaves are {\em algebraic} such a bijection does
exist since the orbits corresponding to a given $w \in W \times W$
have the same dimension.

	In the first section we discuss  the Poisson
structures on $G$ defined by classical $r$-matrices of the form $r = a
- u$ where $a = \sum_{\alpha \in \bRp} e_\alpha
\wedge e_{-\alpha} \in \wedge^2 \g$ and $u \in \wedge^2 \h$. Given such an 
$r$ one constructs a Manin triple of Lie groups $(G\times G, G, G_r)$.
Unlike the one-parameter case (where $u=0$), the dual group $G_r$ will
generally not be an algebraic subgroup of $G\times G$. In fact this
happens if and only if $u \in \wedge^2 \hQ$. Since the quantized
universal enveloping algebra $\Uq$ is a deformation of the algebra of
functions on the algebraic group $G_r$ \cite{DP}, this explains the
difficulty in constructing multi-parameter versions of $\Uq$. From
\cite{LW,STS}, one has that the symplectic leaves are the connected
components of $G \cap G_rxG_r$ where $x\in G$. Since $r$ is
$H$-invariant, the symplectic leaves are permuted by $H$ with the
orbits being contained in Bruhat cells in $G\times G$ indexed by $W
\times W$. In the case where $G_r$ is algebraic, the symplectic leaves
are also algebraic and an explicit formula is given for their
dimension.

	The philosophy of \cite{HL1,HL2} was that, as in the case of
enveloping algebras of algebraic solvable Lie algebras, the primitive
ideals of $\Cq$ should be in bijection with the symplectic leaves of G
(in the case $u=0$). Indeed, since the Lie bracket on $\g_r =
\Lie(G_r)$ is the linearization of the Poisson structure on $G$,
$\Prim \cqp$ should resemble $\Prim U(\g_r)$.  The study of the
multi-parameter versions $\Cqp$ is similar to the case of enveloping
algebras of general solvable Lie algebras. In the general case $\Prim
U(\g_r)$ is in bijection with the co-adjoint orbits in $\g_r^*$ under
the action of the `adjoint algebraic group' of $\g_r$, \cite{Di}. It
is therefore natural that, only in the case where the symplectic
leaves are algebraic, does one expect and obtain a bijection between
the symplectic leaves and the primitive ideals.

	In section 2 we define the notion of an $\bL $-bigraded Hopf
$\K$-algebra, where $\bL$ is an abelian group. When $A$ is finitely
generated such bigradings correspond bijectively to morphisms from the
algebraic group $\bold{L}\spcheck$ to the (algebraic) group $R(A)$ of one-
dimensional representations of $A$. For any antisymmetric bicharacter
$p$ on $\bL$, the multiplication in $A$ may be twisted to give a new
Hopf algebra $A_p$.  Moreover, given a pair of $\bL$-bigraded Hopf
algebras $A$ and $U$ equipped with an $\bL$-compatible Hopf pairing $A
\times U \to \K$, one can deform the pairing to get a new Hopf pairing
between $A_{p^{-1}}$ and $U_p$. This deformation commutes with the
formation of the Drinfeld double in the following sense. Suppose that
$A$ and $U$ are bigraded Hopf algebras equipped with a compatible Hopf
pairing $A^{\text{op}} \times U \to \K$.  Then the Drinfeld double $A\Join U$
inherits a bigrading such that $(A \Join U)_p
\cong A_p \Join
U_p$.

	Let $\Cq$ denote the usual one-parameter quantum group (or
quantum function algebra) and let $\Uq$ be the quantized enveloping
algebra associated to the lattice $\bL$ of weights of $G$. Let $\uqbp$
and $\uqbm$ be the usual sub-Hopf algebras of $\Uq$ corresponding to
the Borel subalgebras $\bp$ and $\bm$ respectively. Let $\Dq = \uqbp
\Join \uqbm$ be the Drinfeld double. Since the groups of  one-dimensional 
representations of $\uqbp$, $\uqbm$, $\Dq$ and $\Cq$ are all
isomorphic to $H= \bold{L}\spcheck$, these algebras are all equipped with
$\bL$-bigradings. Moreover the Rosso-Tanisaki pairing is compatible
with the bigradings on $\uqbp$ and $\uqbm$. For any anti-symmetric
bicharacter $p$ on $\bL$ one may therefore twist simultaneously the
Hopf algebras $\uqbp$, $\uqbm$ and $\Dq$ in such a way that $\Dqpp
\cong U_{q,p}(\frak{b}^+ ) \Join U_{q,p}(\frak{b}^-)$. The algebra
$\Dqpp$ is the `multi-parameter quantized universal enveloping
algebra' constructed by Okado and Yamane \cite{OY} and previously in special 
cases in
\cite{CR,Su}. The canonical
pairing between $\Cq$ and $\Uq$ induces a $\bL$-compatible pairing
between $\Cq$ and $\Dq$. Thus there is an induced pairing between the
multi-parameter quantum group $\Cqp$ and the multi-parameter double
$\Dqp$.  Recall that the Hopf algebra $\Cq$ is defined as the
restricted dual of $\Uq$ with respect to a certain category $\cal{C}$
of modules over $\Uq$. There is a natural deformation functor from
this category to a category $\cal{C}_p$ of modules over $\Dqp$ and
$\Cqp$ turns out to be the restricted dual of $\Dqp$ with respect to
this category.  This Peter-Weyl theorem for $\cqp$ was also found by
Andruskiewitsch and Enriquez in \cite{AE} using a different construction of 
the quantized
universal enveloping algebra and in special cases in \cite{CM,Ha}.

	The main theorem describing the primitive spectrum of $\Cqp$
is proved in the final section. Since $\Cqp$ inherits an
$\bL$-bigrading, there is a natural action of $H$ as automorphisms of
$\Cqp$. For each $w \in W \times W$, we construct an algebra $A_w=
(\Cqp/I_w)_{\cal{E}_w}$ which is a localization of a quotient of
$\Cqp$. For each prime $P \in \Spec \Cqp$ there is a unique $w \in W
\times W$ such that $P \supset I_w$ and $PA_w$ is proper. Thus $\Spec
\cqp \cong \bigsqcup_{w \in W \times W} \Spec _w \Cqp$ where $\Spec _w
\Cqp \cong \Spec A_w$ is the set of primes of type $w$. The key results are then
Theorems \ref{thm 4.14} and \ref{thm 4.15} which state that an ideal
of $A_w$ is generated by its intersection with the center and that $H$
acts transitively on the maximal ideals of the center. From this it
follows that the primitive ideals of $\Cqp$ of type $w$ form an orbit
under the action of $H$.

	An earlier version of our approach to the proof of Joseph's theorem is contained in the unpublished article \cite{HL3}. The approach presented here is a generalization of this proof to the multi-parameter case.

	These results are algebraic analogs of results of Levendorskii
\cite{Le,LS} on the irreducible representations of multi-parameter
function algebras and compact quantum groups. The bijection between
symplectic leaves of the compact Poisson group and irreducible
$*$-representations of the compact quantum group found by Soibelman in
the one-parameter case, breaks down in the multi-parameter case.

	After this work was completed, the authors became aware of the work of Constantini and Varagnolo \cite{CV1,CV2} which has some overlap with the results in this paper.

\bigskip


\section{Poisson Lie Groups} \label{1}

\smallskip

\subsection{Notation} \label{1.1}

Let $\g$ be a complex semi-simple Lie algebra associated to a  Cartan matrix
$[a_{ij}]_{1 \leqslant i,j \leqslant n}$. Let $\{d_i\}_{1 \leqslant
i \leqslant n}$ be relatively prime positive integers such that $[d_i
a_{ij}]_{1 \leqslant i,j \leqslant n}$ is symmetric positive definite.

Let $\h$ be a Cartan subalgebra of $\g$, $\bR$ the associated root
system, $\bB = \{ \alpha_1, \dots, \alpha_n\}$ a basis of $\bR$,
$\bRp$ the set of positive roots and $W$ the Weyl group. 
We denote by $\bP$ and $\bQ$ the
lattices of weights and roots respectively. The fundamental
weights are denoted by $\varpi_1, \dots, \varpi_n $ and the set of
dominant integral weights by $\bP^+ = \sum_{i=1}^n \N \varpi_i$.  
Let $(-,-)$ be  a non-degenerate $\g$-invariant symmetric bilinear form 
on $\g$; it will identify $\g$, resp. $\h$, with its dual $\g^*$, resp.
$\h^*$. The form $(-,-)$ can be chosen in order to induce  a perfect
pairing $\bP \times \bQ \rightarrow \Z$  such that 
$$(\varpi_i,\alpha_j) = \delta_{ij} d_i, \quad (\alpha_i,\alpha_j) =
d_i a_{ij}.$$
Hence $d_i = (\alpha_i,\alpha_i)/2$ and $(\alpha,\alpha) \in 2 \Z$ for
all $\alpha \in \bR$.  
For each $\alpha_j$ we denote by $h_j \in \h$ the corresponding
coroot: $\varpi_i(h_j) = \delta_{ij}$. 
We also set
$$
\npm = \oplus_{\alpha \in \bold{R}_+} \g_{\pm \alpha}, \quad \bpm = \h
\oplus \npm, \quad
\dgoth = \g \times \g, \quad \tgoth = \h \times \h, \quad \upm =
\frak{n}^\pm \times \frak{n}^\mp.
$$

Let $G$ be a connected complex semi-simple algebraic group such
that $\Lie(G) = \g$ and set $D = G \times G$. We  identify $G$ (and
its subgroups) with the diagonal copy inside $D$. We  denote by
$\exp$ the exponential map from $\dgoth$ to $D$. We
shall in general denote a Lie subalgebra of
$\dgoth$ by a gothic symbol and the corresponding connected analytic
subgroup of $D$ by a capital letter.

\smallskip

\subsection{Poisson Lie group structure on $G$} \label{1.2}

Let $a = \sum_{\alpha \in \bRp} e_\alpha \wedge e_{-\alpha} \in \wedge^2
\g$ where the $e_\alpha$ are root vectors such that
$(e_\alpha,e_{\beta}) = \delta_{\alpha,-\beta}.$ Let $u \in \wedge^2
\h$ and set $r = a - u$. Then it is well known that $r$ satisfies the
modified Yang-Baxter equation \cite{BD,Le} and that therefore the tensor
$\pi(g) = (l_g)_* r - (r_g)_* r$ furnishes $G$ with the structure of a
Poisson Lie group, see \cite{D1,LW,STS} ($(l_g)_*$ and $(r_g)_*$ are the
differentials of the left and right translation by $g \in G$). 

We may write $u= \sum_{1 \le i,j \le n} u_{ij} h_i
\otimes h_j$ for a skew-symmetric $n \times n$ matrix
$[u_{ij}]$. The element $u$ can be considered either as an alternating
form on $\h^*$ or a linear map $u \in \End \h$ by the formula
$$\forall \, x \in \h, \quad u(x) = \sum_{i,j} u_{i,j} (x,h_i) h_j.$$
The  Manin triple associated to the Poisson Lie structure on G given by
$r$ is described as follows. Set $u_\pm = u \pm I \in \End \h$ and
define
\begin{gather*}
\vart : \h \hookrightarrow \tgoth, \quad \vart (x) = -(u_-(x),u_+(x)),\\
\agoth = \vart (\h), \quad \g_r = \agoth \oplus \up.
\end{gather*}
Following \cite{STS} one sees easily that the associated Manin triple is
$(\dgoth, \g, \g_r)$ where $\g$ is identified with the diagonal copy
inside $\dgoth.$
Then the corresponding triple of Lie groups is $(D, G, G_r)$,
where $A = \exp(\agoth)$ is an analytic torus and $G_r = A U^+$.
Notice that $\g_r$ is a solvable, but not in general algebraic, Lie
subalgebra of $\frak{d}$.



The following is an easy consequence of the definition of
$\agoth$ and the identities $u_+ + u_- = 2u, u_+ - u_- = 2I$:
\begin{equation}  \label{lem 1.1}
\agoth = \{(x,y) \in \tgoth \; \mid \; x+y = u(y-x)\} = \{
(x,y) \in \tgoth \; \mid \; u_+(x) = u_-(y)\}.
\end{equation}



\smallskip

Recall that $\exp : \h \rightarrow H$ is surjective; let $L_H$ be  its
kernel. 
We shall denote by $\Char{K}$  the group of characters of
an algebraic torus $K$. Any $\chi \in \Char{H}$ is given by $\chi(\exp
x) = \exp d\chi(x)$, $x \in \h$,  where $d\chi \in \h^*$ is the differential of
$\chi.$ 
Then
$$\Char{H} \cong \bL = \pol{L_H} := \{\xi \in \h^* \; \mid \; \xi(L_H)
\subset 2i\pi\Z\}.$$
One can show that $\bL$ has a basis consisting of dominant weights.

Recall that if $\tilde{G}$ is a connected simply connected algebraic
group with Lie algebra $\g$ and maximal torus $\tilde{H}$, we have
\begin{gather*}
L_{\tilde{H}} = \pol{\bP} =  \oplus_{j=1}^n 2i\pi\Z h_j, \quad
\Char{\tilde{H}} \cong \bP, \\
\bQ \subseteq \bL \subseteq \bP, \quad \pi_1(G) = L_H / \pol{\bP} \cong
\bP / \bL.
\end{gather*} 
Notice that $L_H/\pol{\bP}$ is a finite group
and $\exp u(L_H)$ is a  subgroup of $H$. We set
\begin{gather*}
\Gamma_0 = \{ (a,a) \in T \; \mid \; a^2 = 1\}, \quad \Delta = \{
(a,a) \in T \; \mid \; a^2 \in \exp u(L_H) \} ,\\
\Gamma = A \cap H = \{ (a,a) \in T \; \mid \; a =\exp x = \exp y, \;
x+y = u(y-x)\}.
\end{gather*}
It is easily seen that $\Gamma = G \cap G_r$.

\smallskip

\begin{prop} \label{prop 1.2}
We have $\Delta = \Gamma .\Gamma_0$. 
\end{prop}

\begin{pf} We obviously have $\Gamma_0 \subset \Delta$. Let $(\exp h,
\exp h) \in \Gamma, \; h \in \h$. By definition there exist $(x,y)
\in \agoth$, $\ell_1,\ell_2 \in L_H$  such that
$$x = h + \ell_1, \; y = h + \ell_2, \quad y + x = u(y-x).$$  
Hence  $y+x =2h + \ell_1 + \ell_2 = u(\ell_2 - \ell_1)$ and $(\exp h)^2
= \exp 2h =\exp u(\ell_2 - \ell_1).$ This shows $(\exp h, \exp h) \in
\Delta$. Thus $\Gamma.\Gamma_0 \subseteq \Delta$. 

Let $(a,a) \in \Delta$, $ a =\exp h, \; h \in \h$. From $a^2 \in \exp
u(L_H)$ we get $\ell, \ell' \in L_H$ such that $2h = u(\ell') + \ell$.
Set $x = h - \ell/2 - \ell'/2, \; y = h + \ell'/2 - \ell/2$. Then $y+x
= u(y-x)$ and we have $\exp (-\ell/2 - \ell'/2) = \exp (\ell'/2 -
\ell/2) $, since $\ell' \in L_H$. If $ b = \exp(-\ell'/2 + \ell/2)$ we
obtain  $\exp x = \exp y = a b^{-1}$, hence $(a,a) = (\exp x, \exp
y).(b,b) \in \Gamma.\Gamma_0$. Therefore $\Gamma.\Gamma_0 = \Delta$.
\end{pf}

\smallskip

\begin{rem} 
When $u$ is ``generic'' $\Gamma_0$ is not contained in $\Gamma$.  For
example, take $G$ to be $SL(3, \C)$ and $u = \alpha(h_1 \otimes h_2 -
h_2 \otimes h_1)$ with $\alpha \notin \Q$.
\end{rem}

\smallskip

Considered as  a Poisson variety, $G$ decomposes as a disjoint union
of symplectic leaves. Denote by  $\Symp G$ the set of these symplectic
leaves. Since $r$ is $H$-invariant, translation by an element of $H$ is
a Poisson morphism and hence there is an induced action of $H$ on
$\Symp G$.  The key to classifying the symplectic leaves is the
following result, cf. \cite{LW,STS}.

\smallskip

\begin{thm} \label{thm 1.3}
The symplectic leaves of $G$ are exactly the connected components of
$G \cap G_r x G_r$ for $x \in G$. 
\end{thm}

\smallskip

Remark that $A$, $\Gamma$ and $G_r$ are in general
not closed subgroups of $D$. This has for consequence that the
analysis of $\Symp G$ made in \cite[Appendix A]{HL1} in the case $u=0$
does not apply in the general case. 

Set $ Q = H G_r = T U^+$. Then $Q$ is a Borel subgroup of $D$ and,
recalling that the Weyl group associated to the pair $(G,T)$ is $W
\times W$,  the corresponding  Bruhat decomposition yields $D = \sqcup_{w
\in W \times W} Q w Q = \sqcup_{w \in W \times W} Q w G_r.$ Therefore
any symplectic leaf is contained in a Bruhat cell
$QwQ$ for some $w \in W \times W$. 

\smallskip

\begin{defn} A leaf $\cal{A}$ is said to be of type $w$ if $\cal{A}
\subset QwQ$. The set of leaves of type $w$  is denoted by $\Symp_w
G$. 
\end{defn}

\smallskip

For each $w \in W \times W$ set $w= (\wplus,\wminus), w_{\pm} \in W,$
and fix a representative $\wdot$ in the normaliser of $T$.
One shows as in \cite[Appendix A]{HL1} that $G \cap Q \wdot G_r \neq
\emptyset $, for all $w \in W \times W$; hence  $ \Symp_w G \neq
\emptyset $ and $G \cap G_r \wdot G_r \neq
\emptyset$, since $QwQ = \cup_{h \in H} h G_r \wdot G_r.$ 

  The adjoint action of $D$ on itself is denoted by $\Ad$.  Set
\begin{gather*}
U^-_w = \Ad w(U) \cap U^+, \quad A'_w = \{a \in A \;
\mid \; a \wdot G_r = \wdot G_r\}, \\ 
T'_w = \{t \in T \; \mid \; t
G_r \wdot G_r = G_r \wdot G_r\}, \quad H'_w = H \cap T'_w.
\end{gather*}   
Recall that $U^-_w$ is isomorphic to $\C^{l(w)}$ where $l(w)=
l(\wplus) + l(\wminus)$ is the length of $w$. We set $s(w) = \dim
H'_w.$ 

\smallskip

\begin{lem} \label{lem 1.4}
\rom{(i)} $A'_w = \Ad w (A) \cap A$ and  $T'_w = A. \Ad w(A) = A H'_w$.

\rom{(ii)}   We have  $\Lie (A'_w) = \frak{a}'_w = \{ \vart (x)
\, \mid \, x \in \Ker(u_-  w_-^{-1} u_+ - u_+ w_+^{-1} u_-
)\}$ and  $\dim \frak{a}'_w = n - s(w).$ 
\end{lem}

\begin{pf} (i) The first equality is obvious and the second is an easy
consequence of the bijection, induced by 
multiplication, between  $ U^-_w \times T \times U^+ $ and $ QwQ = Q w
G_r .$   

(ii) By definition we have $\frak{a}'_w = \{ \vart (x) \, \mid \, x
\in \h,\, w^{-1}(\vart (x)) \in \agoth \}.$ From  \eqref{lem 1.1}  
we deduce that $\vart (x) \in \frak{a}'_w $ if and only if $u_+
w_+^{-1}(- u_-(x)) = u_- w_-^{-1}( - u_+(x))$.

It follows from (i) that $\dim T'_w = n + \dim H'_w = 2 n - \dim
A'_w$, hence $\dim \frak{a}'_w = n -s(w)$. 
\end{pf}

\smallskip

Recall that $u \in \End \h$ is an alternating bilinear form on 
$\h^*$.
It is easily seen that 
$\forall \, \lambda, \mu \in \h^*,$ $u(\lambda,\mu) = 
-(^tu(\lambda),\mu)$, where $^tu \in \End \h^*$ is the transpose 
of $u$.

\smallskip

\begin{notation} Set
 $^tu = -\Phi, \, \Phi_\pm = \Phi \pm I, \, \sigma(w) = \Phim
w_- \Phip - \Phip \wplus \Phim$,
where $w_\pm \in W$ is considered as an element of $\End \h^*$. 
\end{notation}
\smallskip

Observe that $^tu_{\pm} = -\Phi_\mp$ and that 
\begin{equation} \label{Phi}
u(\lambda,\mu) = (\Phi\lambda,\mu), \quad  \text{for all
$\lambda,\mu \in \h^*$}.
\end{equation} 
 Furthermore, since the
transpose of $w_\pm \in \End \h^*$ is $w_\pm^{-1} \in \End \h$, we have 
$^t\sigma(w) = u_-  w_-^{-1} u_+ - u_+ w_+^{-1} u_-$. Hence by Lemma
\ref{lem 1.4} 
\begin{equation} \label{s(w)}
s(w) = \codim \Ker_{\frak{h}^*} \sigma(w), \quad \dim A'_w = 
\dim \Ker_{\frak{h}^*} \sigma(w).
\end{equation}

\smallskip

\subsection{The algebraic case} \label{1.3}

As explained  in \ref{1.1} the Lie algebra $\g_r$ is in
general not algebraic. We now   describe its algebraic closure.
Recall that a Lie subalgebra $\m$ of $\dgoth$ is said to be algebraic if $\m$
is the Lie algebra of a closed (connected) algebraic subgroup of $D$.

\smallskip

\begin{defn} Let $\m$ be a Lie subalgebra of $\dgoth.$ The smallest
algebraic Lie subalgebra of $\dgoth$ containing $\m$ is called the algebraic
closure of $\m$ and will be denoted by $\tilde{\m}$. 
\end{defn}

\smallskip

Recall that $\g_r = \agoth \oplus \up $. Notice that $\up$ is an algebraic Lie
subalgebra of $\dgoth$; hence it follows from \cite[Corollary II.7.7]{B} that
$\tilde{\g}_r = \atilde \oplus \up.$ Thus we only need to
describe $\atilde$. Since $\tgoth$ is algebraic we have $\atilde
\subseteq \tgoth$ and we are reduced to characterize the algebraic
closure of a Lie subalgebra of $\tgoth = \Lie (T)$.    

The group $T = H \times H$ is an algebraic torus (of rank $2n$).  The
map $\chi \mapsto  d\chi$ identifies $\Char{T}$ with  $\bL
\times \bL$.  

 Let $\k \subset \tgoth$ be a  subalgebra. We set
$$\orth{\k} = \{ \theta \in \Char{T} \; \mid \; \k \subset \Ker_{\frak{t}}
\theta \}.$$
The following  proposition is well known. It can for instance be
deduced  from the results in \cite[ II. 8]{B}.

\smallskip

\begin{prop} \label{prop 1.5}
 Let $\k$ be a  subalgebra of $\tgoth.$ Then
$\tilde{\k} = \cap_{\theta \in \frak{k}^{\perp}}
\Ker_{\frak{t}} \theta$
and $\tilde{\k}$ is the Lie algebra of the closed connected algebraic
subgroup $\tilde{K} = \cap_{\theta \in \frak{k}^{\perp}} \Ker_T \theta.$
\end{prop}

\smallskip

\begin{cor} \label{cor 1.6} 
We have
\begin{gather*}
\frak{a}^\perp = \{ (\lambda,\mu) \in \Char{T} \;
\mid \; \Phip \lambda + \Phim \mu = 0\}, \\
 \tilde{\agoth} = \cap_{(\lambda,\mu) \in \frak{a}^\perp}
\Ker_{\frak{t}} (\lambda,\mu), \quad \tilde{A} = \cap_{(\lambda,\mu)
\in \frak{a}^\perp} \Ker_T (\lambda,\mu).
\end{gather*}
\end{cor} 

\begin{pf} From the definition of $\agoth = \vart (\h)$ we obtain
$$(\lambda,\mu) \in \frak{a}^\perp \Longleftrightarrow
 \forall \, x \in \h, \ \, \lambda(-u_-(x)) + \mu(-u_+(x)) = 0.$$
The first equality then follows from $^tu_\pm = - \Phi_\mp.$ The
remaining assertions are consequences of Proposition \ref{prop 1.5}.
\end{pf}

\smallskip

Set
\begin{gather*}
\hQ =  \Q \otimes_{\Bbb Z} \pol{\bP} = \oplus_{i=1}^n \Q h_i, 
\quad \hQdual = 
\Q \otimes_{\Bbb Z} \bP = \oplus_{i=1}^n \Q \varpi_i \\
\agoth_{\Bbb Q}^\perp = \Q \otimes_{\Bbb Z} \agoth^\perp =\{
(\lambda,\mu) \in \hQdual \times \hQdual \; \mid \; \Phip \lambda +
\Phim \mu = 0\}. 
\end{gather*}
Observe that $\dim_{\Bbb Q} \agoth_{\Bbb Q}^\perp = \rkZ
\agoth^\perp$ and that, by Corollary \ref{cor 1.6},  
\begin{equation} \label{dim agoth}
\dim \tilde{\agoth} = 2n -  \dim_{\Bbb Q} \agoth_{\Bbb Q}^\perp.
\end{equation} 

\smallskip

\begin{lem} \label{lem 1.7}
 $ \agoth_{\Bbb Q}^\perp \cong \{ \nu \in \hQdual \; \mid
\; \Phi\nu \in \hQdual\}.$
\end{lem}

\begin{pf} Define a $\Q$-linear map 
$$\{ \nu \in \hQdual \; \mid \; \Phi\nu \in \hQdual\}
\longrightarrow \agoth_{\Bbb Q}^\perp, \quad \nu \mapsto (-\Phim \nu,
\Phip \nu).$$
It is easily seen that this provides the desired isomorphism.
\end{pf}

\smallskip

\begin{thm} \label{thm 1.8} 
The following assertions are equivalent:

\rom{(i)} $\g_r$ is an algebraic Lie subalgebra of $\dgoth$;


\rom{(ii)} $ u(\bP \times \bP) \subset \Q;$


\rom{(iii)} $\exists \, m \in \N^*, \ \, \Phi(m \bP ) \subset \bP;$


\rom{(iv)} $\Gamma$ is a finite subgroup of $T$.
\end{thm}

\begin{pf} Recall that $\g_r$ is algebraic if and only if $\agoth =
\tilde{\agoth}$, i.e. $n = \dim \agoth = \dim \tilde{\agoth}.$ By
\eqref{dim agoth} and Lemma \ref{lem 1.7}  this is equivalent to
$\Phi(\bP) \subset \hQdual = \Q \otimes_{\Bbb Z} \bP$. The equivalence
of (i) to (iii) then follows from the definitions, \eqref{Phi} and the fact that $^tu
= -\Phi$. 

To prove the equivalence with (iv) we first observe  that, by
Proposition \ref{prop 1.2}, $\Gamma$ is finite if and only if $\exp
u(L_H)$ is finite. Since $L_H/\pol{\bP}$ is finite this is also equivalent to
$\exp u(\pol{\bP})$ being finite. This holds if and only if $u(m
\pol{\bP}) \subset \pol{\bP}$  for some $m \in \N^*$. Hence the result.
\end{pf}

\smallskip

When the equivalent assertions of Theorem \ref{thm 1.8} hold, we shall
say that we are in the {\em algebraic case} or that $u$ {\em is
algebraic}. In this case all the subgroups previously introduced are
closed algebraic subgroups of $D$ and we may define the algebraic
quotient varieties $D/G_r$ and $\Gbar = G/ \Gamma$. Let $p$ be the
projection $G \rightarrow \Gbar$.  Observe that $\Gbar$ is open in in
$D/G_r$ and that the Poisson bracket of $G$ passes to $\Gbar$. We set
\begin{gather*}
\Cwd = G_r \wdot G_r /G_r, \quad \Cw = QwG_r/G_r = \cup_{h \in H}h \Cwd
\\
\Bwd = \Cwd \cap \Gbar , \quad \Bw = \Cw \cap \Gbar, \quad \Aw =
p^{-1}(\Bw).
\end{gather*}

The next theorem  summarizes the description of the symplectic leaves
in the algebraic case.

\smallskip

 \begin{thm} \label{thm 1.9}
\rom{1.} $\Symp_w G \neq \emptyset$ for all $w \in W
\times W$, $\Symp G = \sqcup_{w \in W \times W} \Symp_w G.$ 

\rom{2.} Each symplectic leaf of $\Gbar$, resp. $G$, is of the form $h
\Bwd$, resp. $h \Awd$, for some $h \in H$ and $w \in W \times W$, where
$\Awd$ denotes a fixed connected component of $p^{-1}(\Bwd)$.  

\rom{3.} $\Cwd \cong A_w \times U^-_w$ where $A_w = A/A'_w$ is a torus of rank
$s(w)$. Hence
$\dim \Cwd = \dim \Bwd = \dim \Awd = l(w) + s(w)$ and $H/\Stab_H
\, \Awd$ is a torus of rank $n -s(w).$
\end{thm}

\begin{pf} The proofs are similar to those given in \cite[Appendix
A]{HL1} for the case $u=0$. 
\end{pf}

\medskip


\section{Deformations of Bigraded Hopf Algebras}

\smallskip

\subsection{Bigraded Hopf Algebras and their deformations} \label{2.1}
Let $\bL$ be an (additive) abelian group. We will say that a Hopf
algebra $(A,i,m,\epsilon,\Delta,S)$ over a field $\K$ is an {\em
$\bL$-bigraded Hopf algebra} if it is equipped with an $\bL \times \bL$
grading 
$$A= \bigoplus_{(\lambda,\mu) \in \bold{L} \times \bold{L}}
A_{\lambda,\mu}$$ 
such that
\begin{enumerate}
\item $\K \subset A_{0,0}, \ A_{\lambda,\mu}A_{\lambda',\mu'}
\subset A_{\lambda+\lambda',\mu+\mu'}$ (i.e. $A$ is a graded algebra)
\item $\Delta(A_{\lambda,\mu}) \subset \sum_{\nu \in \bold{L}}
A_{\lambda,\nu}\otimes A_{-\nu,\mu}$ 
\item  $\lambda \neq -\mu$ implies $\epsilon(A_{\lambda,\mu})=0 $
\item $S(A_{\lambda,\mu}) \subset A_{\mu, \lambda}.$
\end{enumerate}

For sake of simplicity we shall often make the following  abuse
of notation: If $a \in A_{\lambda,\mu}$ we will write 
$\Delta(a) = \sum_{\nu} a_{\lambda,\nu} \otimes
a_{-\nu,\mu}$, $a_{\lambda,\nu} \in A_{\lambda,\nu}$, $a_{-\nu,\mu}\in
A_{-\nu,\mu}.$

Let $p : \bL \times \bL \rightarrow \K^*$ be an antisymmetric
bicharacter on $\bL$ in the sense that $p$ is multiplicative in both
entries and that, for all $\lambda,\mu \in \bL$, 
$$\text{(1)} \ p(\mu,\mu) = 1 \ \, ; \ \, \text{(2)} \ p(\lambda,\mu) =
p(\mu,-\lambda).$$ 
Then the map $\ptil : (\bL \times \bL) \times (\bL \times \bL)
\rightarrow \K^*$ given by 
$$\ptil((\lambda,\mu),(\lambda',\mu')) =
p(\lambda,\lambda')p(\mu,\mu')^{-1}$$
is a $2$-cocycle on $\bL \times \bL$ such that $\ptil(0,0) = 1$.

One may then define a new multiplication, $m_p$,  on $A$ by
\begin{equation} \label{m_p}
\forall \,  a \in A_{\lambda,\mu}, \; b \in A_{\lambda',\mu'}, \quad
a\cdot b = p(\lambda,\lambda')p(\mu,\mu')^{-1}ab. 
\end{equation}

\smallskip

\begin{thm} \label{thm 2.1}
$A_p := (A,i,m_p,\epsilon,\Delta,S)$ is an $\bL$-bigraded Hopf algebra.
\end{thm}

\begin{pf}
The proof is a slight generalization of that given in \cite{AST}. It
is well known that $A_p = (A,i,m_p)$ is an associative algebra. Since
$\Delta$ and $\epsilon$ are unchanged, $(A,\Delta,
\epsilon)$ is still a coalgebra. Thus it remains to check that
$\epsilon,\Delta $ are algebra morphisms  and that $S$ is an
antipode. 

Let $x \in A_{\lambda,\mu}$ and $y \in A_{\lambda',\mu'}$. Then
\begin{align*}
\epsilon(x\cdot y) &=  p(\lambda,\lambda')p(\mu,\mu')^{-1}\epsilon(xy) \\
		&=  p(\lambda,\lambda')p(\mu,\mu')^{-1} \delta_{\lambda,-\mu}
\delta_{\lambda',-\mu'} \epsilon(x)\epsilon(y)\\
		&=
p(\lambda,\lambda')p(-\lambda,-\lambda')^{-1}\epsilon(x)\epsilon(y)\\ 
		&= \epsilon(x)\epsilon(y)
\end{align*}
So $\epsilon $ is a homomorphism. Now suppose that $\Delta(x) = \sum
x_{\lambda,\nu}\otimes x_{-\nu,\mu}$ and $\Delta(y) = \sum
y_{\lambda',\nu'}\otimes y_{-\nu',\mu'}$. Then
\begin{align*}
\Delta(x)\cdot \Delta(y) &= (\sum x_{\lambda,\nu}\otimes x_{-\nu,\mu}) \cdot 
(\sum y_{\lambda',\nu'}\otimes y_{-\nu',\mu'})\\
		&= \sum x_{\lambda,\nu}\cdot y_{\lambda',\nu'}\otimes x_{-
\nu,\mu}\cdot y_{-\nu',\mu'}\\
		&=  p(\lambda,\lambda')p(\mu,\mu')^{-1} \sum
p(\nu,\nu')^{-1}p(-\nu,- 
\nu') x_{\lambda,\nu}y_{\lambda',\nu'}\otimes x_{-\nu,\mu}y_{-\nu',\mu'}\\
		&= p(\lambda,\lambda')p(\mu,\mu')^{-1} \Delta(xy)\\
		&= \Delta(x\cdot y)
\end{align*}
So $\Delta$ is also a homomorphism. Finally notice that 
\begin{align*}
\sum S(x_{(1)})\cdot x_{(2)} &= \sum  S(x_{\lambda,\nu}) \cdot x_{-\nu,\mu}\\
		&= \sum  p(\nu,-\nu)p(\lambda,\mu)^{-1} S(x_{\lambda,\nu}) x_{-
\nu,\mu}\\
		&= p(\lambda,\mu)^{-1} \sum  S(x_{\lambda,\nu}) \cdot
x_{-\nu,\mu}\\ 
		&= p(\lambda,\mu)^{-1} \epsilon(x)\\
		&= \epsilon(x)
\end{align*}
A similar calculation shows that $\sum x_{(1)}\cdot S(x_{(2)}) = \epsilon 
(x)$. Hence $S$ is indeed an antipode.
\end{pf}

\smallskip

\begin{rem} The isomorphism class of the algebra $A_p$ depends only on
the cohomology class $[\ptil] \in H^2(\bL \times \bL,\K^*)$, \cite[\S
3]{AST}. 
\end{rem}

\begin{rem}  Theorem \ref{2.1} is a particular case of the following more
general construction. Let $(A,i,m)$ be a $\K$-algebra. Assume that $F
\in GL_{\Bbb K}(A \otimes A)$ is given such that (with the usual notation)

(1) $F  (m \otimes 1) = (m \otimes 1) F_{23} F_{13} \; ; \; F(1
\otimes m) = (1 \otimes m) F_{12} F_{13}$

(2) $F(i \otimes 1) = i \otimes 1 \; ; \; F(1 \otimes i) = 1 \otimes i$

(3) $F_{12} F_{13} F_{23} = F_{23} F_{13} F_{12}$, i.e. $F$ satisfies
the Quantum Yang-Baxter Equation. 

\noindent Set $m_F = m \circ F$. Then
$(A,i,m_F)$ is a $\K$-algebra.

\noindent Assume furthermore that $(A,i,m,\epsilon,\Delta,S)$ is a Hopf algebra
and that

(4) $F : A \otimes A \to A \otimes A $ is morphism of coalgebras

(5) $m F (S \otimes 1) \Delta = m (S \otimes 1) \Delta \; ; \; m F(1
\otimes S) \Delta = m (1 \otimes S)\Delta.$

\noindent Then $A_F := (A,i,m_F,\epsilon,\Delta,S)$ is a Hopf algebra. The
proofs are straightforward verifications and are left to the
interested reader.

When $A$ is an $\bL$-bigraded Hopf algebra and $p$ is an antisymmetric
bicharacter as above, we may define $F \in GL_{\Bbb K}(A \otimes A)$
by
$$\forall a \in A_{\lambda,\mu}, \; \forall b \in A_{\lambda',\mu'},
\; F(a \otimes b) = p(\lambda,\lambda')p(\mu,\mu')^{-1} a \otimes b.$$
It is easily checked that $F$ satisfies the conditions (1) to (5) and
that the Hopf algebras $A_F$ and $A_p$ coincide.

A related construction of the quantization of a monoidal category is
given in \cite{Ly}.
\end{rem}

\smallskip

\subsection{Diagonalizable subgroups of $R(A)$} \label{2.2}
In the case where $\bL$ is a finitely generated group and $A$ is a
finitely generated algebra (which is the case for the multi-parameter
quantum groups considered here), there is a simple geometric
interpretation of $\bL$-bigradings. They correspond to algebraic
group maps from the diagonalizable group $\bL\spcheck$ to the group of
one dimensional representations of $A$.

Assume that $\K$ is algebraically closed.
Let $(A,i,m,\epsilon,\Delta,S)$ be a Hopf $\K$-alge\-bra. Denote by $R(A)$
the multiplicative gro\-up of one dimensional repre\-sen\-tations of
$A$, i.e. the cha\-racter group of the algebra $A$.  Notice that when
$A$ is a finitely generated $\K$-algebra, $R(A)$ has the structure of an
af\-fine algebraic group over $\K$, with algebra of regular functions
given by $\K[R(A)] = A/J$ where $J$ is the semi-prime ideal $\cap_{h
\in R(A)} \Ker h$. Recall that there are two natural group homomorphisms $l, r :
R(A) \rightarrow
\Aut_{\Bbb K}(A) $ given by
\begin{gather*}
l_h(x) = \sum h(S(x_{(1)}))x_{(2)} = \sum h^{-1}(x_{(1)})x_{(2)} \\
r_h(x) = \sum x_{(1)}h(x_{(2)}).
\end{gather*}

\smallskip

\begin{thm} Let $A$ be a finitely generated Hopf algebra and let $\bL$ be a 
finitely
generated abelian group. Then there is a natural bijection between:
\begin{enumerate}
\item $\bL$-bigradings on $A$;
\item Hopf algebra maps $A \to \K \bL$ (where $\K \bL$ denotes the group 
algebra);
\item morphisms of algebraic groups $\bL\spcheck \to R(A)$.
\end{enumerate}
\end{thm}

{\renewcommand{\qedsymbol}{}
\begin{pf} The bijection of the last two sets of maps is well-known. Given 
an $\bL$-bigrading on $A$, we may define a map $\phi: A \to \K \bL$ by 
$\phi(a_{\lambda,\mu}) =
\epsilon(a)u_\lambda$. It is easily verified that this is a Hopf algebra 
map. Conversely,
given a map $\bL\spcheck \to R(A)$ we may construct an $\bL$ bigrading using the
following result.
\end{pf}}

\smallskip

\begin{thm} \label{thm 2.2}
Let $(A,i,m,\epsilon,\Delta,S)$ be a finitely generated Hopf algebra
over $\K$. Let $H$ be a closed diagonalizable algebraic subgroup of
$R(A)$. Denote by $\bL$ the (additive) group of characters of $H$ and by
$\pair{-}{-} : \bL \times H \rightarrow \K^*$ the natural pairing. For
$(\lambda,\mu) \in \bL \times \bL$ set
$$ A_{\lambda,\mu} = \{ x \in A \, \mid  \, \forall h \in H, \; l_h(x) = 
\pair{\lambda}{h} x, \; r_h(x) = \pair{\mu}{h} x \}. $$
Then  $(A,i,m,\epsilon,\Delta,S)$ is an $\bL$-bigraded Hopf algebra.
\end{thm}

\begin{pf}
Recall that any element of $A$ is contained in a finite dimensional subcoalgebra
 of $A$.  Therefore the actions of $H$ via $r$ and $l$
are locally finite. Since they commute and $H$ is diagonalizable, $A$
is $\bL \times \bL$ diagonalizable. Thus the decomposition $A=
\bigoplus_{(\lambda,\mu) \in \bold{L} \times \bold{L}}A_{\lambda,\mu}$
is a grading.

Now let $C$ be a finite dimensional subcoalgebra of $A$ and let $\{
c_1, \dots , c_n\} $ be a basis of $H \times H$ weight vectors.
Suppose that $\Delta(c_i) = \sum t_{ij} \otimes c_j$. Then since $c_i
= \sum t_{ij}\epsilon(c_j)$, the $t_{ij}$ span $C$ and it is easily
checked that $\Delta(t_{ij}) = \sum t_{ik} \otimes t_{kj}$.  Since
$l_h(c_i) = \sum h^{-1}(t_{ij})c_j$ for all $h \in H$ and the $c_i$
are weight vectors, we must have that $h(t_{ij}) = 0$ for $i \neq j$.
This implies that $$l_h(t_{ij}) = h^{-1}(t_{ii})t_{ij}, \quad
r_h(t_{ij}) = h(t_{jj})t_{ij}$$ and that the map $\lambda_i(h) =
h(t_{ii})$ is a character of $H$.  Thus $t_{ij} \in
A_{-\lambda_i,\lambda_j}$ and hence $$\Delta(t_{ij}) = \sum
t_{ik}\otimes t_{kj} \in \sum A_{-\lambda_i,\lambda_k}\otimes
A_{-\lambda_k,\lambda_j}.$$ This gives the required condition on
$\Delta$.  If $\lambda + \mu \neq 0$ then there exists an $ h \in H$
such that $\pair{-\lambda}{h} \neq \pair{\mu}{h}$. Let $x \in
A_{\lambda,\mu}$. Then $$ \pair{\mu}{h} \epsilon(x) = \epsilon(r_h(x))
= h(x) =
\epsilon(l_{h^{-1}}(x)) = \pair{-\lambda}{h} \epsilon(x).$$
Hence $\epsilon(x) =0$. The  assertion on $S$  follows similarly.
\end{pf} 

\smallskip

\begin{rem} In particular, if $G$ is any algebraic group and $H$ is a 
diagonalizable subgroup with character group $\bL$, then we may deform
the Hopf algebra $\K[G]$ using an antisymmetric bicharacter on $\bL$.
Such deformations are algebraic  analogs of the deformations studied
by Rieffel in \cite{Re}. 
\end{rem}

\smallskip

\subsection{Deformations of dual pairs} \label{2.3}
Let $A$ and $U$ be a dual pair of Hopf algebras. That is, there exists a 
bilinear pairing $\langle\; \mid \; \rangle : A \times U \to \K$ such that:
\begin{enumerate}
\item $\hpair{a}{1} = \epsilon(a) \; ; \; \hpair{1}{u} = \epsilon(u)$
\item $ \langle a \mid u_1u_2 \rangle = \sum \langle a_{(1)} \mid u_1 \rangle 
\langle a_{(2)} \mid u_2 \rangle$
\item $\langle a_1a_2 \mid u\rangle = \sum \langle a_1 \mid
u_{(1)}\rangle
 \langle a_2\mid u_{(2)} \rangle$
\item $\langle S(a) \mid u \rangle = \langle a\mid S(u) \rangle .$
\end{enumerate}

Assume that $A$ is bigraded by $\bL$, $U$ is bigraded by an abelian
group $\bQ$ and that there 
is a homomorphism $ \, \breve{}: \bQ \to \bL$ such that 
\begin{equation} \label{pairing condition}
\langle A_{\lambda, \mu}\mid U_{\gamma,\delta}\rangle \neq 0 \ \text{
only if } \  \lambda +\mu = \breve{\gamma}+\breve{\delta}.
\end{equation}
In this case we will call the pair $\{A,U\}$ an {\em $\bL$-bigraded
dual pair.} We shall be interested in \S 3 and \S 4 in the case where
$\bQ = \bL$ and $\, \breve{} \, = Id$. 

\smallskip

\begin{rem} Suppose that  the bigradings  above are induced from subgroups 
$H$ and $\breve{H}$ of $R(A)$ and $R(U)$ respectively and that the map $\bQ \to
\bL$ is induced from a map $h \mapsto \breve{h}$ from $H$ to $\breve{H}$.
Then the condition on
the pairing may be restated as the fact that the form is ad-invariant
in the sense that for all $a \in A$, $u \in U$ and $h \in H$,
$$\langle \ad_{h} a\mid  u \rangle = \langle a\mid
\ad_{\breve{h}}u \rangle$$ 
where $\ad_h a = r_h l_h(a)$.
\end{rem}

\smallskip

\begin{thm} \label{thm 2.3}
Let $\{A,U\}$ be the bigraded dual pair as described above.  Let $p$
be an antisymmetric bicharacter on $\bL$ and let $\breve{p}$ be the
induced bicharacter on $\bQ$.  Define a bilinear form
$\langle\;\mid\;\rangle_p:A_{p^{-1}}\times U_{\breve{p}} \to \K$ by: 
$$\langle a_{\lambda,\mu}\mid u_{\gamma,\delta}\rangle_p =
p(\lambda,\breve{\gamma})^{-1} p(\mu,\breve{\delta})^{-1} \langle
a_{\lambda,\mu}\mid u_{\gamma,\delta}\rangle.$$ 
Then $\langle
\;\mid\;\rangle_p$ is a Hopf pairing and
$\{A_{p^{-1}},U_{\breve{p}}\}$
 is an $\bL$-bigraded dual pair.
\end{thm}

\begin{pf}
Let $a \in A_{\lambda, \mu}$ and let $u_i \in U_{\gamma_i,\delta_i}$,
$i=1,2.$ Then 
$$ \langle a\mid u_1u_2 \rangle_p =
p(\breve{\gamma}_1,\breve{\gamma}_2)p(\breve{\delta}_1,\breve{\delta}_2)^{-1}
p(\lambda,\breve{\gamma}_1+\breve{\gamma}_2)^{-
1}p(\mu,\breve{\delta}_1+\breve{\delta}_2)^{-1} \langle a\mid u_1u_2
\rangle.$$ 
Suppose that $\Delta(a) = \sum_\nu a_{\lambda,\nu}\otimes
a_{-\nu,\mu}$.  Then by the assumption on the pairing, the only
possible value of $\nu$ for which $\langle a_{\lambda,\nu}\mid u_1
\rangle \langle a_{-\nu,\mu}\mid u_2 \rangle$ is non-zero is $\nu =
\breve{\gamma}_1 + \breve{\delta}_1 - \lambda = \mu - \breve{\gamma}_2 -
\breve{\delta}_2$. Therefore
\begin{multline*} \langle a_{(1)}\mid u_1 \rangle_p \langle
a_{(2)}\mid  u_2 \rangle_p  = 
p(\lambda,\breve{\gamma}_1)^{-1} p(\nu,\breve{\delta}_1)^{-1}p(-
\nu,\breve{\gamma}_2)^{-1} p(\mu,\breve{\delta}_2)^{-1} 
	\langle a_{(1)}\mid u_1 \rangle \langle a_{(2)}\mid u_2 \rangle \\
= p(\lambda,\breve{\gamma}_1)^{-1} p(\mu-\breve{\gamma}_2-
\breve{\delta}_2,\breve{\delta}_1)^{-1}p(\lambda-\breve{\gamma}_1-
\breve{\delta}_1,\breve{\gamma}_2)^{-1} p(\mu,\breve{\delta}_2)^{-1} 
	\langle a_{(1)}\mid u_1 \rangle \langle a_{(2)}\mid u_2 \rangle\\
= p(\breve{\gamma}_1,\breve{\gamma}_2)p(\breve{\delta}_1,\breve{\delta}_2)^{-1}
p(\lambda,\breve{\gamma}_1+\breve{\gamma}_2)^{-
1}p(\mu,\breve{\delta}_1+\breve{\delta}_2)^{-1}
\langle a\mid u_1u_2 \rangle = \langle a\mid u_1u_2 \rangle_p.
\end{multline*}
This proves the first axiom. The others are verified similarly.
\end{pf}

\smallskip

\begin{cor} \label{cor 2.4} 
Let $\{ A,U,p \}$ be as in Theorem \ref{thm 2.3}. Let $M$ be a right
$A$-comodule with structure map $\rho : M \to M \otimes A$.  Then $M$
is naturally endowed with $U$ and $U_{\breve{p}}$ left module
structures, denoted by $(u,x) \mapsto ux$ and $(u,x) \mapsto u \cdot
x$ respectively. Assume that $M = \oplus_{\lambda \in \bold{L}}
M_\lambda$ for some
$\K$-subspaces such that $\rho(M_\lambda) \subset \sum_\nu M_{-\nu} \otimes
A_{\nu,\lambda}$. Then  we have $U_{\gamma,\delta} M_{\lambda} 
\subset M_{\lambda - \breve{\gamma} - \breve{\delta}}$ and 
the two structures are related by  $$\forall u \in U_{\gamma,\delta}, \; \forall x \in M_\lambda, \quad
u\cdot x = p(\lambda,\breve{\gamma} - \breve{\delta})
p(\breve{\gamma},\breve{\delta}) ux.$$
\end{cor}
 
\begin{pf} Notice that the coalgebras $A$ and $A_{p^{-1}}$ are the same. Set
$\rho(x) = \sum x_{(0)} \otimes x_{(1)}$ for all $x \in 
M$. Then it is easily checked that the following formulas define the
desired $U$ and $U_{\breve{p}}$ module structures: 
$$\forall u \in U, \quad ux =
\sum x_{(0)} \hpair{x_{(1)}}{u}, \quad u \cdot x =
\sum x_{(0)} \hpair{x_{(1)}}{u}_p.$$
When $x \in M_\lambda$ and $u \in U_{\gamma,\delta}$ the additional
condition yields 
$$u \cdot x = \sum x_{(0)} p(\nu,-\breve{\gamma})p(\lambda,-\breve{\delta})
\hpair{x_{(1)}}{u}.$$
 But $\hpair{x_{(1)}}{u} \ne 0$  forces $-\nu =
\lambda -\breve{\gamma} - \breve{\delta}$, hence
$u \cdot x =  p(\lambda,\breve{\gamma}-
\breve{\delta})p(\breve{\gamma},\breve{\delta}) \sum  x_{(0)} 
\hpair{x_{(1)}}{u} = p(\lambda,\breve{\gamma}-
\breve{\delta})p(\breve{\gamma},\breve{\delta}) ux.$ 
\end{pf}

\smallskip

Denote by $\Aop$ the opposite algebra of the $\K$-algebra $A$. Let
$\{\Aop,U, \hpair{ \, }{ \, }\}$ be a dual pair of Hopf algebras. The double
$A \Join U$ is defined as follows, \cite[3.3]{DCL}. Let $I$ be the
ideal of the tensor algebra $T(A \otimes U)$ generated by elements of
type
\begin{gather*}
1 - 1_A, \quad 1 - 1_U \tag{a} \\
x x'  - x \otimes x', \ x,x' \in A, \quad    y y' - y
\otimes y', \ y,y' \in U \tag{b} \\
x_{(1)} \otimes y_{(1)} \hpair{x_{(2)}}{y_{(2)}} - \hpair{x_{(1)}}{y_{(1)}}
y_{(2)} \otimes  x_{(2)}, \ x \in A, \, y \in U \tag{c}
\end{gather*} 
Then the algebra $A \Join U : = T(A \otimes U)/I$ is  called {\em the
 Drinfeld double of $\{A,U\}$}. It is a Hopf algebra in a natural way:
\begin{gather*}
\Delta(a \otimes u) = (a_{(1)} \otimes u_{(1)}) \otimes (a_{(2)}
\otimes u_{(2)}), \\ \epsilon(a \otimes u) = \epsilon(a)\epsilon(u), \quad
S(a \otimes u) = (S(a) \otimes 1)(1 \otimes S(u)).
\end{gather*}

Notice for further use that $A \Join U$ can equally be defined by
relations of type (a), (b), (c$_{x,y}$)  or (a), (b), (c$_{y,x}$),
where we set
\begin{align*}
x \otimes y  &=  \hpair{x_{(1)}}{y_{(1)}} \hpair{x_{(3)}}{S(y_{(3)})} y_{(2)}
 \otimes x_{(2)}, \ x \in A, \, y \in U \tag{c$_{x,y}$} \\
y\otimes x &=  \hpair{x_{(1)}}{S(y_{(1)})} \hpair{x_{(3)}}{y_{(3)}} x_{(2)}
\otimes y_{(2)}, \ x \in A, \, y \in U \tag{c$_{y,x}$}
\end{align*}

\smallskip

\begin{thm} \label{thm 2.5}
Let $\{ \Aop,U \}$ be an $\bL$-bigraded
dual pair, $p$ be an antisymmetric bicharacter on $\bL$ and 
$\breve{p}$ be the induced bicharacter on $\bQ$. Then $A \Join U$
inherits an $\bL$-bigrading and there is a natural
isomorphism of $\bL$-bigraded Hopf algebras:
$$(A \Join U)_p \cong A_p \Join U_{\breve{p}}.$$
\end{thm}

\begin{pf} Recall that as a $\K$-vector space $A \Join U$ identifies
with $A \otimes U$. Define an $\bL$-bigrading on $A \Join U$ by
$$\forall \alpha, \beta \in \bL, \quad (A \Join U)_{\alpha,\beta} =
\sum_{\lambda-\breve{\gamma} = \alpha, \mu - \breve{\delta} = \beta}
A_{\lambda,\mu}  \otimes U_{\gamma, \delta}.$$
To verify that this yields a structure of graded algebra on $A \Join
U$ it suffices to check that the defining relations of $A \Join U$ are
homogeneous. This is clear for relations of type (a) or (b). Let
$x_{\lambda,\mu} \in A _{\lambda,\mu}$ and $y_{\gamma,\delta} \in
U_{\gamma,\delta}.$ Then the corresponding relation of type (c) becomes
\begin{equation*}
\sum_{\nu,\xi} x_{\lambda,\nu} y_{\gamma,\xi}
\hpair{x_{-\nu,\mu}}{y_{-\xi,\delta}} -
\hpair{x_{\lambda,\mu}}{y_{\gamma,\xi}}  y_{-\xi,\delta} x_{-\nu,\mu}.
\tag{$\star$}
\end{equation*}
When a term of this sum is non-zero we obtain
$-\nu + \mu = -\breve{\xi} + \breve{\delta},$ $\lambda + \nu =
\breve{\gamma} + \breve{\xi}.$
Hence $\lambda - \breve{\gamma} = - \nu + \breve{\xi} = -\mu +
\breve{\delta}$, which shows that the relation ($\star$) is
homogeneous. It is easy to see that the conditions (2), (3), (4) of
\ref{2.1} hold. Hence $A \Join U$ is an $\bL$-bigraded Hopf algebra. 

Notice that $(A_p)^{\text{op}} \cong (\Aop)_{p^{-1}}$, so that Theorem
\ref{thm 2.3} defines a suitable pairing between $(A_p)^{\text{op}}$ and
$U_{\breve{p}}$. Thus $A_p \Join U_{\breve{p}}$ is defined. Let $\phi$
be the natural surjective homomorphism  from $T(A \otimes U)$ onto
$A_p \Join U_{\breve{p}}$. 
To check that $\phi$ induces an isomorphism it again suffices to check that
 $\phi$ vanishes on the defining relations of $(A \Join U)_p$. Again, this 
is easy for relations of type (a) and (b). The relation
($\star$) says that 
\begin{equation*}
p(\lambda,\breve{\gamma})p(-\nu,\breve{\xi})\hpair{x_{-\nu,\mu}}{y_{-\xi,\delta}}
x_{\lambda,\nu}\cdot y_{\gamma,\xi} - p(\breve{\xi},\nu)
p(\breve{\delta}, -\mu) \hpair{x_{\lambda,\mu}}{y_{\gamma,\xi}}
y_{-\xi,\delta} \cdot x_{-\nu,\mu} =0  
\end{equation*}
in $(A \Join U)_p$.
Multiply the left hand side of this equation by $p(\lambda,-\breve{\gamma})
p(\mu,-\breve{\delta})$ and apply $\phi$. We obtain the following
expression in $A_p \Join U_{\breve{p}}$:
\begin{equation*}
p(-\nu,\breve{\xi}) p(\mu,-\breve{\delta}) \hpair{x_{-\nu,\mu}}{y_{-\xi,\delta}}
x_{\lambda,\nu} y_{\gamma,\xi} - p(\lambda, -\breve{\gamma})
p(\nu,-\breve{\xi}) \hpair{x_{\lambda,\mu}}{y_{\gamma,\xi}}
y_{-\xi,\delta}  x_{-\nu,\mu} 
\end{equation*}
which is equal to
$$\hpair{x_{-\nu,\mu}}{y_{-\xi,\delta}}_p x_{\lambda,\nu}
y_{\gamma,\xi} - \hpair{x_{\lambda,\mu}}{y_{\gamma,\xi}}_p
y_{-\xi,\delta}  x_{-\nu,\mu}.$$
But this is a defining relation of type (c) in $A_p \Join
U_{\breve{p}}$, hence zero.

It remains to see that $\phi$
induces an isomorphism of Hopf algebras, which is a straightforward
consequence of the definitions.   
\end{pf}

\smallskip

\subsection{Cocycles} \label{2.4}
Let  $\bL$ be, in this section, an arbitrary  free abelian group with basis
$\{\omega_1, \dots, \omega_n\}$ and set $\h^* = \C \otimes_{\Bbb Z}
\bL.$  We  freely use the terminology  of
\cite{AST}. Recall that  $H^2(\bL,\C^*)$ is in bijection with the set 
$\cal{H}$ of multiplicatively antisymmetric $n \times n$-matrices
$\gamma=[\gamma_{ij}]$. This bijection maps the class $[c]$ onto the matrix
defined by $\gamma_{ij} =
c(\omega_i,\omega_j)/c(\omega_j,\omega_i)$. Furthermore it is
an isomorphism of groups with respect to component-wise multiplication
of matrices.

\smallskip

\begin{rem} The notation is as in \ref{2.1}.  We recalled that  the
isomorphism class of the algebra $A_p$ depends only on the cohomology
class $[\ptil] \in H^2(\bL \times \bL,\K^*)$. Let $\gamma \in \cal{H}$
be the matrix associated to $p$ and $\gamma^{-1}$ its inverse in
$\cal{H}$.   Notice that the
multiplicative matrix associated to  $[\ptil]$ is then 
$\tilde{\gamma} = \smatrix{\gamma}{1}{1}{\gamma^{-1}}$ 
in the basis given by the $(\omega_i,0), (0,\omega_i) \in \bL
\times \bL$. Therefore the isomorphism class of the algebra $A_p$
depends only on the cohomology class $[p] \in H^2(\bL,\K^*)$.
\end{rem} 

\smallskip

Let $\hbar \in \C^*$. If $x \in \C$ we set $\qexp{x} = \exp(-x
\hbar/2)$. In particular $q = \exp(-\hbar/2)$.
Let $ u : \bL \times \bL \rightarrow \C$ be a complex alternating
$\Z$-bilinear form.
Define
\begin{equation} \label{defn of p}
p: \bL \times \bL \rightarrow \C^*, \quad p(\lambda,\mu) =
\exp \big(-\frac{\hbar}{4} u(\lambda,\mu)\big)
=\qexp{\frac{1}{2}u(\lambda,\mu)}.
\end{equation}
Then it is clear that $p$ is an antisymmetric  bicharacter on $\bL$.

Observe that, since $\frak{h}^* = \C \otimes_{\Bbb Z} \bL$, there is a
natural isomorphism of additive groups between $\wedge^2 \h$ and the
group of complex alternating $\Z$-bilinear forms on $\bL$, where $\h$ is
the $\C$-dual of $\frak{h}^*$.  Set 
$\cal{Z}_{\hbar} =
\{ u \in \wedge^2 \h \; \mid \; u(\bL \times \bL) \subset 
\frac{ 4 i \pi}{\hbar} \Z \}.$

\smallskip

\begin{thm} \label{thm 2.6}
There are isomorphisms of abelian groups:
$$H^2(\bL,\C^*) \cong \cal{H} \cong \wedge^2 \h /\cal{Z}_\hbar.$$
\end{thm}

\begin{pf}  The first isomorphism has been described above. 
Let $\gamma =[\gamma_{ij}] \in \cal{H}$ and choose $u_{ij}$, $1 \le i
< j \le n$ such that $\gamma_{ij} = \exp( - \frac{\hbar}{2} u_{ij})$.
We can define $u \in \wedge^2 \h$ by setting $u(\omega_i,\omega_j) =
u_{ij}$, $1 \le i < j \le n$.
It is then easily seen that one can  define an injective 
 morphism of abelian groups 
$$\varphi :
H^2(\bL,\C^*) \cong
\cal{H} \longrightarrow \wedge^2 \h /
\cal{Z}_\hbar, \quad  \varphi(\gamma) = [u]$$
where $[u]$ is the class of $u$.  If $u \in \wedge^2 \h $, define a
$2$-cocycle $p$ by the formula \eqref{defn of p}. Then the
multiplicative matrix associated to $[p] \in H^2(\bL,\C^*)$ is given
by 
$$\gamma_{ij} = p(\omega_i,\omega_j)/p(\omega_j,\omega_i) =
p(\omega_i,\omega_j)^2 =
\exp(- \frac{\hbar}{2} u(\omega_i,\omega_j)).$$
This shows that $[u] = \varphi([\gamma_{ij}])$; thus $\varphi$ is an
isomorphism.
\end{pf}

\smallskip
 
We list some consequences of Theorem \ref{thm 2.6}. We denote by $[u]$
an element of $\wedge^2 \h / \cal{Z}_\hbar$ and we set $[p] =
\varphi^{-1}([u])$. We have seen that we can define  a representative
$p$ by the formula \eqref{defn of p}.

1. $[p] \ \text{ of finite order in $H^2(\bL,\C^*)$} \Leftrightarrow
u(\bL \times \bL) \subset \frac{i \pi}{\hbar} \Q$, and 
$q \ \text{root of unity} \Leftrightarrow \hbar \in i \pi \Q.$

2. Notice that  $u=0$ is algebraic, whether $q$ is a root of unity or
not.
 Assume that $q$ is a root of unity; then we get from 1 that
$$[p] \ \text{ of finite order} \Leftrightarrow  u \ \text{is
algebraic}.$$ 

3. Assume that $q$ is not a root of unity and that $u \ne 0$. Then
$[p]$ of finite order implies $(0) \ne u(\bL \times \bL) \subset
\frac{i \pi}{\hbar}\Q$. This shows that
$$0 \ne u \ \text{algebraic} \Rightarrow [p] \ \text{ is not of finite
order}.$$

\smallskip

\begin{defn} 
The bicharacter $p: (\lambda,\mu) \mapsto q^{\frac{1}{2}u(\lambda,\mu)}$ is
called $q$-rational if $u \in \wedge^2 \h$ is algebraic.
\end{defn}

The following technical result will be used in the next section. Recall the definition of 
$\Phi_-=\Phi-I$ given in the Section 1.

\begin{prop} \label{Kprop}
Let $ \bK = \left\{\lambda \in \bL \, : \, (\Phi_-\lambda,\bL) \subset
\frac{4i\pi}{\hbar}\Z \right\}$. If $q$ is not a root of unity, then $\bK  = 0$.
\end{prop}

\begin{pf}
Let $\lambda \in \bK$. We can define $z : \h^*_\Q \to \Q$, by
$$
\forall \, \mu \in \h^*_\Q, \quad (\Phi_-\lambda,\mu) = \cte 
z(\mu).  
$$
The map $z$ is clearly $\Q$-linear. It follows, since
$(\phantom{x},\phantom{y})$ is non-degenerate on $\h^*_\Q$, that there
exists $\nu \in \h^*_\Q$ such that $z(\mu) = (\nu,\mu)$ for
all $\mu \in \h^*_\Q$. 
Therefore  $\Phi_- \lambda = \cte \nu$, and so $\Phi \lambda = \lambda +
\cte \nu$.

Now, $(\Phi \lambda,\lambda) = u(\lambda,\lambda) = 0$ implies that 
$$
\cte (\nu,\lambda) = -(\lambda,\lambda)
$$
If $(\lambda,\lambda) \neq 0$ then $\hbar \in i\pi \Q$, 
contradicting the assumption that $q$ is not a root of unity. Hence
$(\lambda,\lambda) = 0$, which forces $\lambda = 0$ since $\lambda \in
\bL \subset \h^*_\Q$.
\end{pf}

\medskip


\section{Multiparameter Quantum Groups} 

\smallskip

\subsection{One-parameter quantized enveloping algebras} \label{3.1}
The notation is as in sections 1 and 2. In particular we fix a lattice
$\bL$ such that $\bQ \subset \bL \subset \bP$ and we denote by $G$ the
connected semi-simple algebraic group with maximal torus $H$  such
that $\Lie(G) = \g$ and $\Char{H} \cong \bL$. 

Let $q \in \C^*$ and assume that {\em $q$ is not a root of unity.} Let
$\hbar \in \C \setminus i \pi \Q$ such that $q = \exp(-\hbar/2)$ as in
\ref{2.4}. We set
$$q_i = q^{d_i}, \quad \hat{q_i} = (q_i - q_i^{-1})^{-1}, \quad 1 \le
i \le n.$$
Denote by $U^0$ the group algebra of  $\Char{H}$, hence
$$U^0 = \C[k_\lambda \, ; \, \lambda \in \bL], \quad k_0 =1, \quad
k_\lambda k_\mu = k_{\lambda + \mu}.$$
Set $k_i = k_{\alpha_i}$, $ 1 \le i \le n$. The one parameter quantized
enveloping algebra associated to this data, cf. \cite{T}, is the Hopf
algebra 
$$U_q(\g) = U^0[e_i, f_i \, ; \, 1 \le i \le n]$$
with defining  relations:
\begin{gather*}
k_\lambda e_j k_\lambda^{-1} = q^{(\lambda,\alpha_j)} e_j, \quad
k_\lambda f_j k_\lambda^{-1} = q^{-(\lambda,\alpha_j)} f_j \\
 e_if_j - f_je_i = \delta_{ij} \hat{q_i} (k_i-k_i^{-1}) \\
\sum_{k = 0}^{1-a_{ij}} (-1)^{k}
\left[ \begin{smallmatrix} 1-a_{ij} \\ k \end{smallmatrix} \right]_{q_i}
 e_i^{1-a_{ij}-k}e_je_i^k = 0, \mbox{ if } i \neq j \\
\sum_{k = 0}^{1-a_{ij}} (-1)^{k}
\left[ \begin{smallmatrix} 1-a_{ij} \\ k \end{smallmatrix} \right]_{q_i}
 f_i^{1-a_{ij}-k}f_jf_i^k = 0,
\mbox{ if } i \neq j
\end{gather*}
where $[m]_t =(t-t^{-1}) \ldots (t^m-t^{-m})$ and
$\left[ \begin{smallmatrix} m \\ k \end{smallmatrix} \right]_t = 
\frac{[m]_t}{[k]_t [m-k]_t}.$ The Hopf algebra structure is given by
\begin{gather*}
\Delta(k_\lambda) = k_\lambda \otimes k_\lambda, \quad
\epsilon(k_\lambda) = 1, \quad S(k_\lambda) = k_\lambda^{-1} \\
\Delta(e_i) = e_i \otimes 1 + k_i \otimes e_i, \quad 
\Delta(f_i) = f_i \otimes k_i^{-1} + 1 \otimes f_i \\
\epsilon(e_i) = \epsilon(f_i) =0, \quad S(e_i) = -k_i^{-1} e_i, \quad 
S(f_i) = -f_i k_i.
\end{gather*}
We define subalgebras of $\Uq$ as follows
\begin{gather*}
 \uqnp = \C[e_i, \, ; \, 1 \le i \le n], \quad \uqnm = \C[f_i, \, ; \, 1
\le i \le n] \\
\uqbp = U^0[e_i, \, ; \, 1 \le i \le n], \quad \uqbm =  U^0[f_i, \, ;
\, 1 \le i \le n].
\end{gather*}
For simplicity we shall set $U^{\pm} =  U_q(\frak{n}^{\pm})$. Notice that
$U^0$ and $U_q(\frak{b}^{\pm})$ are Hopf subalgebras of $\Uq.$ Recall \cite{L}
that the multiplication in $\Uq$ induces isomorphisms of vector spaces
$$\Uq \cong U^- \otimes U^0 \otimes U^+ \cong U^+ \otimes U^0 \otimes
U^-.$$
Set  $\bQ_+ = \oplus_{i=1}^n \N \alpha_i$ and 
$$ \forall \beta \in \bQ_+, \quad  U^{\pm}_\beta = \{ u \in U^\pm \,
\mid \, \forall \lambda \in \bL, \,  k_\lambda u k_\lambda^{-1} =
q^{(\lambda,\pm\beta)} u \}.$$
Then one gets: $U^\pm = \oplus_{\beta \in \bold{Q}_+} U^\pm_{\pm
\beta}$.

\smallskip

\subsection{The Rosso-Tanisaki-Killing form} \label{3.2}
Recall the following result, \cite{Ro1,T}.

\smallskip

\begin{thm} \label{thm 3.1}
\rom{1.} There exists a unique non degenerate Hopf pairing 
$$\hpair{\,}{\,} : \uqbp^{\text{op}} \otimes \uqbm \longrightarrow
\C$$
satisfying the following conditions:

\rom{(i)} $\hpair{k_\lambda}{k_\mu} = q^{-(\lambda,\mu)};$

\rom{(ii)} $ \forall\lambda \in \bL, 1 \le i \le n,$
$\hpair{k_\lambda}{f_i} = \hpair{e_i}{k_\lambda} = 0;$

\rom{(iii)} $ \forall 1 \le i,j \le n$, $\hpair{e_i}{f_j} = -
\delta_{ij} \hat{q_i}.$

\rom{2.} If $\gamma, \eta \in \bQ_+$,
$\hpair{U^+_\gamma}{U^-_{-\eta}} \neq 0$ implies $ \gamma = \eta$. 
\end{thm}

\smallskip

The results of \S \ref{2.3} then apply and we may define the
associated double:
$$D_q(\g) = \uqbp \Join \uqbm.$$
It is well known, e.g.  \cite{DCL}, that
$$\Dq = \C[s_\lambda, t_\lambda, e_i, f_i \, ; \, \lambda \in \bL, 1
\le i \le n]$$
where $s_\lambda = k_\lambda \otimes 1$, $t_\lambda = 1 \otimes
k_\lambda$, $e_i = e_i \otimes 1$, $f_i = 1 \otimes f_i$. The defining
relations of the double given in \S \ref{2.3} imply that
\begin{gather*}
s_\lambda t_\mu = t_\mu s_\lambda, \quad e_i f_j - f_j e_i =
\delta_{ij} \hat{q_i}(s_{\alpha_i} - t_{\alpha_i}^{-1}) \\
s_\lambda e_j s_\lambda^{-1} = q^{(\lambda,\alpha_j)} e_j, \ t_\lambda
e_j t_\lambda^{-1} = q^{(\lambda,\alpha_j)} e_j, \ s_\lambda f_j
s_\lambda^{-1} = q^{-(\lambda,\alpha_j)} f_j, \  t_\lambda f_j
t_\lambda^{-1} = q^{-(\lambda,\alpha_j)} f_j. 
\end{gather*}
 It follows that
$$\Dq/(s_\lambda - t_\lambda \, ; \,  \lambda \in \bL) \isomto \Uq, \
e_i \mapsto e_i, \, f_i 
\mapsto f_i,\, s_\lambda \mapsto k_\lambda, \, t_\lambda \mapsto
k_\lambda.$$
Observe that this yields  an isomorphism of Hopf algebras.   
The next proposition collects some  well known elementary facts.

\smallskip

\begin{prop}  \label{prop 3.2}
\rom{1.} Any finite dimensional simple $U_q(\frak{b}^\pm)$-module is
one dimensional and $R(U_q(\frak{b}^\pm))$ identifies with $H$ via
$$\forall h \in H, \quad h(k_\lambda) = \pair{\lambda}{h}, \quad h(e_i) =
0, \quad   h(f_i) = 0.$$

\rom{2.} $R(\Dq)$ identifies with $H$ via
$$\forall h \in H,  \quad h(s_\lambda) = \pair{\lambda}{h}, \quad  h(t_\lambda) =
\pair{\lambda}{h}^{-1} , \quad h(e_i) = h(f_i) =0.$$
\end{prop}

\smallskip

\begin{cor} \label{cor 3.3}
\rom{1.} $\{\uqbp^{\text{op}}, \uqbm\}$ is an $\bL$-bigraded dual pair.
We have
$$k_\lambda \in U_q(\frak{b}^\pm)_{-\lambda,\lambda}, \quad e_i \in
\uqbp_{-\alpha_i,0}, \quad f_i \in \uqbm_{0,-\alpha_i}.$$

\rom{2.} $\Dq$ is an $\bL$-bigraded Hopf algebra where
$$s_\lambda \in \Dq_{-\lambda,\lambda}, \quad t_\lambda \in
\Dq_{\lambda, - \lambda}, \quad e_i \in \Dq_{-\alpha_i,0}, \quad f_i
\in \Dq_{0,\alpha_i}.$$
\end{cor}

\begin{pf}  1. Observe that for all $h \in H$,
\begin{gather*}
l_h(k_\lambda) = h^{-1}(k_\lambda) = \pair{-\lambda}{h} k_\lambda,
\quad r_h(k_\lambda) = h(k_\lambda) =\pair{\lambda}{h} k_\lambda, \\
l_h(e_i) =  h^{-1}(k_i) e_i = \pair{-\alpha_i}{h} e_i , \quad r_h(e_i)
= e_i, \\
l_h(f_i) = f_i, \quad r_h(f_i) = h(k_i^{-1}) f_i = \pair{-\alpha_i}{h}
f_i.
\end{gather*}
It is then clear that $U^+_{-\gamma,0} = U^+_\gamma$ and
$U^-_{0,-\gamma} = U^-_{-\gamma}$ for all $\gamma \in \bQ_+$.
The claims then follow from these formulas, Theorem \ref{thm 2.2},
Theorem \ref{thm 3.1},  and  the definitions.

2. The fact that $\Dq$ is an $\bL$-bigraded Hopf algebra follows from
Theorem \ref{thm 2.2}. The assertions about the $\bL \times \bL$
degree of the generators is proved by direct computation using
Proposition \ref{prop 3.2}. \end{pf}

\smallskip

\begin{rem} We have shown in Theorem \ref{thm 2.5} that, as a double,
$\Dq$ inherits an $\bL$-bigrading given by:
$$\Dq_{\alpha,\beta} = \sum_{\lambda- \gamma = \alpha, \mu - \delta =
\beta} \uqbp_{\lambda,\mu} \otimes \uqbm_{\gamma,\delta}.$$
It is easily checked that this bigrading coincides with the bigrading
obtained in the above corollary by means of Theorem \ref{thm 2.2}.
\end{rem}

\smallskip

\subsection{One-parameter quantized function algebras} \label{3.3}
 Let $M$ be a left $\Dq$-mo\-du\-le. The dual $M^*$ will be considered
in the usual way
as a left $\Dq$-module by the rule: $(uf)(x)= f(S(u)x)$, $x \in M, f \in
M^*$, $u \in \Dq$. Assume that $M$ is an $\Uq$-module. An element $x
\in M$ is said to have {\em weight $\mu \in
\bL$} if $k_\lambda x = q^{(\lambda,\mu)}x$ for all $\lambda \in \bL$;
 we denote by $M_\mu$ the subspace of elements of weight $\mu$.

It is known, \cite{D1},  that the category of finite dimensional (left)
$\Uq$-modules is a completely reducible braided rigid monoidal category.
Set $\bL^+ = \bL \cap \bP^+$ and recall that for each $\Lambda \in
\bL^+$ there exists a finite dimensional simple module of highest
weight $\Lambda$, denoted by $L(\Lambda)$, cf. \cite{Ro2} for
instance. One has $L(\Lambda)^* \cong L(w_0 \Lambda)$ where $w_0$ is
the longest element of $W$. (Notice that the results quoted usually cover the case where
$\bL=\bQ$. One defines the modules $L(\lambda)$ in the general case in the following way. Let us denote temporarily the algebra $\Uq$ for a given choice of $\bL$ by $U_{q,\bL}(\g)$. Given a module $L(\lambda)$ defined on $U_{q,\bQ}(\g)$ we may define an action of $U_{q,\bL}(\g)$ by setting
$k_\lambda.x = q^{(\lambda,\mu)}x$ for all elements $x$ of weight $\mu$, where $ q^{(\lambda,\mu)}$ is as defined  in section 2.4.)

  Let $\calCq$ be the subcategory of finite dimensional
$\Uq$-modules consisting of finite direct sums of $L(\Lambda)$, $\Lambda
\in \bL^+$. The category  $\calCq$ is closed under tensor products and the
formation of duals. Notice that $\calCq$ can be considered as a
braided rigid monoidal category of $\Dq$-modules where $s_\lambda, t_\lambda$
act as $k_\lambda$ on an object of $\calCq.$

Let $M \in \text{obj}(\calCq),$ then $M
=\oplus_{\mu \in \bold{L}} M_\mu$. For $f \in M^*$, $v \in M$ we define the
 coordinate function  $c_{f,v} \in \Uq^*$ by
$$\forall u \in \Uq, \quad c_{f,v}(u) = \pair{f}{uv}$$
where $\pair{\,}{\,}$ is the duality pairing. Using the standard isomorphism $(M\otimes N)^* \cong N^* \otimes M^*$ one has the following formula for multiplication,
$$ c_{f,v} c_{f',v'} = c_{f' \otimes f, v \otimes v'}.$$

\begin{defn} The quantized function algebra $\Cq$ is the restricted
dual of $\calCq$: that is to say
$$\Cq = \C[c_{f,v} \, ; \, v \in M, f \in M^*, \, M \in
\text{obj}(\calCq)].$$ 
\end{defn}

The algebra $\Cq$ is a Hopf algebra; we denote by $\Delta, \epsilon, S$
the comultiplication, counit and antipode on $\Cq$. If $\{v_1, \dots,
v_s; f_1, \dots, f_s\}$ is a dual basis for $M \in \text{obj}(\calCq)$
one has

\begin{equation} \label{formulas 3.1}
\Delta(c_{f,v}) = \sum_i c_{f,v_i} \otimes c_{f_i,v}, \quad
\epsilon(c_{f,v}) = \pair{f}{v}, \quad S(c_{f,v}) = c_{v,f}.
\end{equation}
Notice that we may assume that $v_j \in M_{\nu_j}, \, f_j \in
M^*_{-\nu_j}$. We set 
$$C(M) = \C \langle c_{f,v} \, ; \, f\in M^* , \, v \in M \rangle,
\quad C(M)_{\lambda,\mu} = \C \langle c_{f,v} \, ; \, f\in
M^*_\lambda, \, v \in M_\mu \rangle.$$
Then $C(M)$ is a  subcoalgebra  of $\Cq$ such that 
$C(M) = \bigoplus_{(\lambda,\mu) \in \bold{L} \times
\bold{L}} C(M)_{\lambda,\mu}.$ When $M = L(\Lambda)$ we abbreviate the
notation to $C(M) = C(\Lambda)$. It is then classical that 
$$\Cq =\bigoplus_{\Lambda \in \bold{L}^+} C(\Lambda).$$
Since $\Cq \subset \Uq^*$ we have a duality pairing 
$$\langle \ , \ \rangle : \Cq \times \Dq \longrightarrow \C.$$ 

Observe that there is a  natural injective morphism of algebraic groups 
$$H \longrightarrow R(\Cq), \quad h(c_{f,v}) =
\pair{\mu}{h}\epsilon(c_{f,v}) \ \text{for all $v \in M_\mu$, $ M \in
\text{obj}(\calCq)$}.$$  
The associated automorphisms $r_h, l_h \in \Aut(\Cq)$ are then
described by
$$\forall c_{f,v} \in C(M)_{\lambda,\mu}, \quad r_h(c_{f,v} ) =
\pair{\mu}{h} c_{f,v}, \ l_h(c_{f,v} ) = \pair{\lambda}{h}c_{f,v}.$$
Define 
$$\forall (\lambda,\mu) \in \bL \times \bL, \quad \Cq_{\lambda,\mu} =
\{a \in \Cq \, \mid \, r_h(a) =\pair{\mu}{h} a, \, l_h(a) =
\pair{\lambda}{h}a \}.$$

\smallskip

\begin{thm}  \label{thm 3.4}
The pair of Hopf algebras $\{\Cq, \Dq\}$ is an
$\bL$-bigraded dual pair.
\end{thm}

\begin{pf} It follows from \eqref{formulas 3.1} that $\Cq$ is an
$\bL$-bigraded Hopf algebra. The axioms (1) to (4) of \ref{2.3} are
satisfied by definition of the Hopf algebra $\Cq$. We take   $\,
\breve{\,} \, $ to be the identity map of $\bL$. 
The condition \eqref{pairing condition} is consequence of 
$\Dq_{\gamma,\delta} M_\mu
\subset M_{\mu - \gamma - \delta}$ for all $M \in \calCq$.
To verify this inclusion, notice that
$$e_j \in \Dq_{-\alpha_j,0}, \ f_j \in \Dq_{0,\alpha_j},
\quad e_j M_\mu \subset M_{\mu + \alpha_j }, \ f_j M_\mu \subset
M_{\mu - \alpha_j }.$$
The result then follows easily 
\end{pf}

\smallskip

Consider the algebras $D_{q^{-1}}(\g)$ and $\C_{q^{-1}}[G]$ and use  
$ \, \hat{} \,$  to distinguish elements, subalgebras, etc. of
$D_{q^{-1}}(\g)$ and $\C_{q^{-1}}[G]$. It is easily verified that the
map $\sigma : \Dq \to D_{q^{-1}}(\g)$ given by
$$s_\lambda \mapsto \hat{s}_\lambda, \ t_\lambda \mapsto
\hat{t}_\lambda, \ e_i \mapsto q_i^{1/2} \hat{f}_i \hat{t}_{\alpha_i},
\ f_i \mapsto q_i^{1/2} \hat{e}_i \hat{s}^{-1}_{\alpha_i}$$
is an isomorphism of Hopf algebras.

For each $\Lambda \in \bL^+$, $\sigma$ gives a bijection $\sigma : L(-w_0 \Lambda) \to \hat{L}(\Lambda)$ which sends $v \in L(-w_0
\Lambda)_\mu $ onto $\hat{v} \in \hat{L}(\Lambda)_{-\mu}.$
Therefore we obtain  an isomorphism $\sigma: \C_{q^{-1}}[G] \to
\Cq$ such that 
\begin{equation} \label{3.3.1}
\forall \, f \in L(-w_0 \Lambda)^*_{-\lambda}, \, v \in L(-w_0
\Lambda)_\mu , \quad \sigma(\hat{c}_{\hat{f}, \hat{v}}) = c_{f,v}.
\end{equation} 
Notice that
\begin{equation} \label{3.3.2}
\sigma(\Dq_{\gamma,\delta}) = D_{q^{-1}}(\g)_{-\gamma, - \delta} \ \
\text{and} \ \
\sigma(\C_{q^{-1}}[G]_{\lambda,\mu} ) = \C_q[G]_{-\lambda,-\mu}.
\end{equation}

\smallskip

\subsection{Deformation of one-parameter quantum groups} \label{3.4}
We continue with the same notation. Let $[p] \in H^2(\bL, \C^*).$ As
seen in \S \ref{2.4} we can, and we do, choose $p$ to be an
antisymmetric bicharacter such that 
$$\forall \lambda, \mu \in \bL,
\quad p(\lambda,\mu) = q^{\frac{1}{2}u(\lambda,\mu)}$$
 for some $u \in \wedge^2 \h$. Recall that $\ptil \in Z^2(\bL
\times\bL, \C^*)$, cf. \ref{2.1}.  

We now apply the results of \S \ref{2.1} to $\Dq$ and
$\Cq$. Using Theorem \ref{thm 2.1} we can  twist  $\Dq$ by
$\ptil^{-1}$ and $\Cq$ by $\ptil$. 
The resulting $\bL$-bigraded Hopf algebras will be denoted by
$\Dqp$ and $\Cqp$. The algebra $\Cqp$ will be referred to as the 
{\em multi-parameter quantized
function algebra}. Versions of $\Dqp$ are referred to by some authors as the
 {\em multi-parameter quantized enveloping algebra}. Alternatively, this name can be applied to the quotient of $\Dqp$ by the radical of the pairing with  $\Cqp$.

\smallskip

\begin{thm}
Let $\uqpbp$ and $\uqpbm$ be the deformations by $p^{-1}$ of $\uqbp$
and $\uqbm$ respectively. Then the deformed pairing
$$\hpair{\,}{\,}_{p^{-1}} : \uqpbp^{\text{op}} \otimes \uqpbm
\rightarrow \C$$ 
is a non-degenerate Hopf pairing satisfying:
\begin{equation} \label{pairing}
\forall \, x \in U^+, \,  y  \in U^-, \, \lambda,\mu \in \bL, \quad
\hpair{x \cdot k_\lambda}{y \cdot k_\mu}_{p^{-1}} = \qexp{(\Phim
\lambda,\mu)} \hpair{x}{y}.
\end{equation}
Moreover,
$$ \uqpbp \Join \uqpbm  \cong (\uqbp \Join \uqbm)_{p^{-1}}= \Dqp  .$$
\end{thm}

\begin{pf}  
By Theorem \ref{thm 2.3} the deformed pairing is given by
$$\hpair{a_{\lambda,\mu}}{u_{\gamma,\delta}}_{p^{-1}} =
p(\lambda,\gamma)p(\mu,\delta)
\hpair{a_{\lambda,\mu}}{u_{\gamma,\delta}}.$$
To prove \eqref{pairing} we can assume that $x \in U^+_{-\gamma,0}, \,
y \in U^-_{0,-\nu}$. Then we obtain
\begin{align*}
\hpair{x \cdot k_\lambda}{y \cdot k_\mu}_{p^{-1}} &= p(\lambda +
\gamma,\mu)p(\lambda,\mu - \nu) \hpair{x \cdot k_\lambda}{y \cdot
k_\mu} \\
&= p(\lambda,2\mu)p(\lambda - \mu,\gamma - \nu)\qexp{-(\lambda,\mu)}
\hpair{x}{y}
\end{align*}
by the definition of the product $ \, \cdot \, $ and \cite[2.1.3]{T}.
But $\hpair{x}{y} = 0$ unless $\gamma = \nu$, hence the result.
Observe in particular that $\hpair{x}{y}_{p^{-1}} =
\hpair{x}{y}$. Therefore \cite[2.1.4]{T} shows that $\hpair{ \, }{ \,
}_{p^{-1}}$ is non-degenerate on  $U^+_\gamma \times U^-_{-\gamma}.$
It then follows  from \eqref{pairing} and Proposition \ref{Kprop} that
$\hpair{\,}{\,}_{p^{-1}}$ is non-degenerate.
The remaining isomorphism follows from Theorem \ref{thm 2.5}. 
\end{pf}  

Many authors have defined multi-parameter quantized enveloping
algebras. In \cite{Ha,OY} a definition is given using
explicit generators and relations, and in \cite{AE} the construction is
made by twisting the comultiplication, following \cite{R}.
 It can be easily verified that these
algebras and the algebras $\Dqp$  coincide. The
construction of a multi-parameter quantized function algebra by
twisting the multiplication was first performed in the $GL(n)$-case in
\cite{AST}. 

The fact that $D_{q,p^{-1}}(\g)$ and $\C_{q,p}[G]$ form a Hopf dual
pair has already been observed in particular cases, see e.g.
\cite{Ha}. We will now deduce from the previous results that this
phenomenon holds for an arbitrary semi-simple group.

\smallskip

\begin{thm} \label{thm 3.5} 
$\{\Cqp, D_{q,p^{-1}}(\g)\}$ is an $\bL$-bigraded dual pair. 
The associated pairing is given by 
$$\forall a \in \Cqp_{\lambda,\mu}, \, \forall u \in \Dqp_{\gamma,\delta},
\quad \pair{a}{u}_p = p(\lambda,\gamma)p(\mu,\delta)\pair{a}{u}.$$
\end{thm}

\begin{pf} This follows from Theorem \ref{thm 2.3} applied to the pair
$\{A,U\} = \{\Cq,\Dq\}$ and the bicharacter $p^{-1}$ (recall that the
map $\, \breve{\,} \, $ is the identity). 
\end{pf}

\smallskip

 Let $M \in \text{obj}(\calCq)$. The left $\Dq$-module structure on
$M$ yields a right $\Cq$-comodule structure in the usual way. Let
$\{v_1, \dots, v_s;f_1, \dots, f_s\}$ be a dual basis for $M$.  The
structure map $\rho : M \to M \otimes \Cq$, is given by $ \rho(x) =
\sum_j v_j \otimes c_{f_j,x}$ for  $x \in M$. 
Using this comodule structure on $M$, one can check  that
$$M_\mu = \{ x \in M \, \mid \, \forall h \in H, \; r_h(x) =
\pair{\mu}{h} x\}.$$ 

\smallskip

\begin{prop} \label{prop 3.6} 
Let  $M \in \text{obj}(\calCq)$. Then
$M$ has a natural structure of left $D_{q,p^{-1}}(\g)$ module. Denote by
$M\spbreve$ this module and by $(u,x) \mapsto u \cdot x$ the action of
$D_{q,p^{-1}}(\g)$. Then
$$\forall u \in \Dq_{\gamma,\delta}, \, \forall x \in M_\lambda, \quad
u \cdot x = p(\lambda,\delta - \gamma)p(\delta,\gamma)ux.$$
\end{prop}

\begin{pf} The proposition is a translation in this particular setting of
Corollary \ref{cor 2.4}.
\end{pf}

\smallskip

Denote by $\calCqp$ the subcategory of finite dimensional left
$D_{q,p^{-1}}(\g)$-modules whose  objects are the $M\spbreve$, $M \in
\text{obj}(\calCq)$. It follows from Proposition \ref{prop 3.6}
that if $M \in\text{obj}(\calCq)$,
then $M\spbreve = \oplus_{\mu \in \bold{L}} M\spbreve_\mu$, where
$$M\spbreve_\mu = \{ x \in M \, \mid \, \forall \alpha \in \bL, \
s_\alpha\cdot x = p(\mu,2\alpha) q^{(\mu,\alpha)}x, \; t_\alpha \cdot x
= p(\mu, -2\alpha) q^{(\mu,\alpha)} x \}.$$
Notice that $p(\mu,\pm 2\alpha) q^{(\mu,\alpha)} = q^{\pm(\Phipm
\mu,\alpha)}.$ 

\smallskip

\begin{thm} \label{thm 3.7}
\rom{1.} The functor $M \to M\spbreve$ from $\calCq$ to
$\calCqp$ is an equivalence of rigid monoidal categories.

\rom{2.} The Hopf pairing $\pair{\,}{\,}_p$ identifies the Hopf
algebra $\Cqp$ with the  restricted dual of $\calCqp$,
i.e. the Hopf algebra of coordinate functions on the objects of
$\calCqp$.
\end{thm}

\begin{pf} 1. One needs in particular to prove that, for $M, N \in
\text{obj}(\calCq)$,  there are natural isomorphisms
of $D_{q,p^{-1}}(\g)$-modules: $\varphi_{M,N} : (M \otimes N)\spbreve
\to M\spbreve \otimes N\spbreve$.  These isomorphisms are given by
$x \otimes y \mapsto p(\lambda,\mu) x
\otimes y$ for all $x \in M_\lambda, \, y \in N_\mu$. The other
verifications are elementary.

2. We have to show that if $M \in \text{obj}(\calCq)$, $f \in M^*$, $v
\in M$ and $u \in D_{q,p^{-1}}(\g)$, then $\pair{c_{f,v}}{u}_p =
\pair{f}{u \cdot v}.$  It suffices to prove the result in the case 
where $f \in M^*_\lambda, \, v \in M_\mu$ and $ u \in
D_{q,p^{-1}}(\g)_{\gamma,\delta}.$ Then
\begin{align*}
\pair{f}{u \cdot v} &= p(\mu, \delta - \gamma)p(\delta,\gamma)
\pair{f}{uv} \\
& = \delta_{-\lambda + \gamma + \delta, \mu} \, p(-\lambda + \gamma +
\delta,\delta - \gamma)p(\delta,\gamma) \pair{f}{uv} \\
&= p(\lambda,\gamma)p(\mu,\delta) \pair{f}{uv} \\ 
&= \pair{c_{f,v}}{u}_p
\end{align*}
by Theorem \ref{thm 3.5}. 
\end{pf}

\smallskip

Recall that we introduced in \S \ref{3.3} isomorphisms $\sigma : \Dq
\to D_{q^{-1}}(\g)$ and $\sigma :  \Cq \to \C_{q^{-1}}[G]$. 
From \eqref{3.3.2} it follows that, after twisting by $\ptil^{-1}$ or
$\ptil$, $\sigma$ induces isomorphisms 
$$D_{q,p^{-1}}(\g) \isomto D_{q^{-1},p^{-1}}(\g), \quad \C_{q^{-1},p}[G]
\isomto \C_{q,p}[G]$$
which satisfy \eqref{3.3.1}.

\smallskip

\subsection{Braiding isomorphisms} \label{3.5} 
We remarked above that the categories $\calCqp$ are braided. In the
one parameter case this braiding is well-known. Let $M$ and $N$ be
objects of $\calCq$.  Let $E: M \otimes N \to M \otimes N$ be the
operator given by 
$$ E(m \otimes n ) = q^{(\lambda,\mu)} m \otimes n$$
for $m \in M_\lambda$ and $n \in N_\mu$. Let $\tau : M \otimes N \to N
\otimes M$ be
the usual twist operator. Finally let $C$ be the operator given by left 
multiplication by
$$ C = \sum_{\beta \in \bQ_+} C_\beta $$
where $C_\beta$ is the canonical element of $\Dq$ associated to the 
non-degenerate
pairing $U^+_\beta \otimes U^-_{-\beta}\to \C $ described above. Then
one deduces from \cite[4.3]{T} that the 
operators
$$\theta_{M,N}= \tau \circ C \circ E^{-1} : M \otimes N \to N \otimes M$$
define the braiding on $\calCq$.

As mentioned above, the category $ \calCqp$ inherits a braiding given
by $$\psi_{M,N} = \varphi_{N,M} \circ \theta_{M,N} \circ
\varphi_{M,N}^{-1}$$
where  $\varphi_{M,N}$ is the isomorphism $ (M \otimes N)\spbreve
\isomto M\spbreve \otimes N\spbreve$ introduced in the proof of
Theorem \ref{thm 3.7} (the same formula can be found in \cite[\S
10]{AE} and in a more general situation in \cite{Ly}).  We now note
that these general operators are of the same form as those in the one
parameter case. Let $M$ and $N$ be objects of $\calCqp$ and let $E: M
\otimes N \to M \otimes N$ be the operator given by
$$ E(m \otimes n ) = q^{(\Phi_+\lambda,\mu)} m \otimes n$$ for $m \in
M_\lambda$ and $n \in N_\mu$.  Denote by $C_\beta$ the canonical
element of $\Dqp$ associated to the nondegenerate pairing
$\uqpbp_{-\beta,0} \otimes \uqpbm_{0,-\beta} \to \C$ and let $C :M
\otimes N \to
M \otimes N$ be the operator given by left multiplication by
$$ C = \sum_{\beta \in \bQ_+} C_\beta .$$

\begin{thm}
The braiding operators $\psi_{M,N}$ are given by
$$\psi_{M,N} = \tau \circ C \circ E^{-1}.$$
Moreover $(\psi_{M,N})^* = \psi_{M^*,N^*}$.
\end{thm}

\begin{pf} The  assertions follow easily from the analogous assertions 
for $\theta_{M,N}$.
\end{pf}

 The following commutation relations are well known \cite{So},
\cite[4.2.2]{LS}. We include a proof for completeness.

\smallskip

\begin{cor} \label{cor 3.10}
Let $\Lambda, \Lambda' \in \bold{L}^+$,  let $g \in
L(\Lambda')^*_{-\eta}$ and $f \in
L(\Lambda)^*_{-\mu}$ and let $v_\Lambda \in L(\Lambda)_\Lambda$. Then for any $v \in
L(\Lambda')_\gamma$,
 $$\cfv{g}{v} \cdot \cfv{f}{v_\Lambda} = q^{(\Phip
\Lambda,\gamma)-(\Phip \mu,\eta)}
\cfv{f}{v_\Lambda} \cdot \cfv{g}{v} + 
q^{(\Phip \Lambda,\gamma) - (\Phip \mu,\eta)} \sum_{\nu \in
\bold{Q}_+}\cfv{f_\nu}{v_\Lambda} \cdot \cfv{g_\nu}{v}$$ 
where $f_\nu \in (\uqpbp f)_{-\mu+\nu}$ and
$g_\nu \in (\uqpbm g)_{-\eta-\nu}$ are such that
$\sum f_\nu \otimes g_\nu = \sum_{\beta \in \bQ^+\backslash \{0\}} C_\beta 
(f\otimes g) $.
\end{cor}

\begin{pf}
Let $\psi = \psi_{L(\Lambda), L(\Lambda')}$. Notice that
$$c_{f\otimes g,\psi(v_\Lambda\otimes v)} = c_{\psi^*(f \otimes g), v_\Lambda \otimes
v}.$$
Using the theorem above we obtain
$$\psi^*(f \otimes g) = q^{-(\Phi_+\mu,\eta)} (g \otimes f + \sum g_\nu 
\otimes f_\nu)$$
and 
\begin{equation} \label{braiding}
\psi(v_\Lambda \otimes v) = q^{-(\Phi_+\Lambda,\gamma)}(v \otimes 
v_\Lambda).
\end{equation}
Combining these formulae yields the required relations.
\end{pf}

\medskip


\section{Prime and Primitive  Spectrum of $\Cqp$}

\smallskip

In this section we prove our main result on the primitive spectrum
of $\Cqp$; namely that the $H$ orbits inside $\Prim_w \Cqp$ are
parameterized by the double Weyl group. For completeness we have
attempted to make the proof more or less self-contained. The overall
structure of the proof is similar to that used in \cite{HL2} except
that the proof of the key \ref{thm 4.12} (and the lemmas leading up to
it) form a modified and abbreviated version of Joseph's proof of this
result in the one-parameter case \cite{J1}. One of the main differences
with the approach of \cite{J1} is the use of the Rosso-Tanisaki form
introduced in \ref{3.2} which simplifies the analysis of the adjoint
action of $\Cqp$. The ideas behind the first few results of this
section go back to Soibelman's work in the one-parameter `compact'
case \cite{So}. These ideas were adapted to the multi-parameter case
by Levendorskii \cite{Le}.

\smallskip

\subsection{Parameterization of the prime spectrum}
Let $q,p$ be as in \S \ref{3.4}. For simplicity we set
$$A = \Cqp$$
and the product $a \cdot b$ as defined in \eqref{m_p} will be
denoted by $ab$. 

For each $\Lambda \in \bL^+$ choose weight vectors 
$$v_\Lambda \in L(\Lambda)_\Lambda, \ v_{w_0 \Lambda} \in
L(\Lambda)_{w_0\Lambda}, \ f_{-\Lambda} \in L(\Lambda)^*_{-\Lambda}, \
f_{-w_0 \Lambda} \in L(\Lambda)^*_{-w_0 \Lambda}$$
such that $\pair{f_{-\Lambda}}{v_\Lambda} = 
\pair{f_{-w_0\Lambda}}{v_{w_0\Lambda}} = 1 $. Set
$$A^+ = \sum_{\mu \in \bL^+} \sum_{f \in L(\mu)^*} \C c_{f,v_\mu},
\quad A^- = \sum_{\mu \in \bL^+} \sum_{f \in L(\mu)^*} \C
c_{f,v_{w_0\mu}}.$$ 
Recall the following result.

\smallskip

\begin{thm} \label{thm 4.1}
 The multiplication map $A^+ \otimes A^- \to A $ is
surjective.
\end{thm}

\begin{pf}  Clearly  it is enough to prove the theorem in the
one-parameter case. 
When $\bL = \bP$   the result
is proved in  \cite[3.1]{So} and \cite[Theorem 3.7]{J1}.

The general case can be deduced from the simply-connected case as
follows.
One first observes that $\Cq \subset \C_q[\tilde{G}] =
\bigoplus_{\Lambda \in \bP^+} C(\Lambda)$. Therefore any $a \in \Cq $
can be written in the form $a = \sum_{\Lambda', \Lambda'' \in \bP^+}
\cfv{f}{v_{\Lambda'}} \cfv{g}{v_{-\Lambda''}}$ where $\Lambda' -
\Lambda'' \in \bL$. Let $\Lambda \in \bP$ and $\{v_i;f_i\}_i$ be a
dual basis of $L(\Lambda)$. Then we have 
$$1 = \epsilon(\cfv{v_\Lambda}{f_{-\Lambda}}) = \sum_i
\cfv{f_i}{v_\Lambda} \cfv{v_i}{f_{-\Lambda}}.$$
Let $\Lambda'$  be as above and choose $\Lambda$ such that $\Lambda +
\Lambda' \in \bL^+$. Then, for all
$i$, $\cfv{f}{v_{\Lambda'}}\cfv{f_i}{v_\Lambda} \in C(\Lambda +
\Lambda') \cap A^+$ and
$\cfv{v_i}{f_{-\Lambda}}\cfv{g}{v_{-\Lambda''}} \in C(-w_0(\Lambda +
\Lambda'')) \cap A^-$. The result then follows by inserting $1$
between the terms $ \cfv{f}{v_{\Lambda'}}$ and 
$\cfv{g}{v_{-\Lambda''}}$. 
\end{pf}

\smallskip

\begin{rem} 
The algebra $A$ is a Noetherian domain (this result will not be used
in the sequel). 
The fact that $A$ is a domain follows from the same
result in \cite[Lemma 3.1]{J1}. The fact that $A$ is Noetherian is
consequence of \cite[Proposition 4.1]{J1} and \cite[Theorem 3.7]{CQ}. 
\end{rem}

\smallskip

For each $y \in W$ define the following ideals of $A$
\begin{gather*}
I_y^+ = \langle \cfv{f}{v_\Lambda} \, \mid f \in (\uqpbp
L(\Lambda)_{y \Lambda})^\perp, \, \Lambda \in \bL^+ \rangle, \\ 
I_y^- = \langle \cfv{f}{v_{w_0\Lambda}} \, \mid f \in (\uqpbm
L(\Lambda)_{y w_0\Lambda} )^\perp, \, \Lambda \in \bL^+
\rangle
\end{gather*}
where $( \ )^\perp$ denotes the orthogonal in $L(\Lambda)^*$. Notice
that $I_y^- =\sigma(\hat{I}^+_y)$, $\sigma$ as in \S \ref{3.4}, and
that $I^\pm_y$ is an $\bL \times \bL$ homogeneous ideal of $A$.

\smallskip

\begin{notation} For $ w = (w_+,w_-) \in W \times W$ set
$I_w = I^+_{w_+} + I^-_{w_-}.$
For $\Lambda \in \bL^+$ set
$c_{w\Lambda} = \cfv{f_{-w_+\Lambda}}{v_\Lambda} \in
C(\Lambda)_{-w_+\Lambda, \Lambda}$  and
$\tilde{c}_{w\Lambda} = \cfv{v_{w_-\Lambda}}{f_{-\Lambda}} \in
C(-w_0\Lambda)_{w_-\Lambda,-\Lambda} $.
\end{notation}

\smallskip

\begin{lem} \label{lem 4.2}
Let $\Lambda \in \bL^+$ and $a \in A_{-\eta,\gamma}$. Then 
\begin{gather*} 
\cw{\Lambda} a \equiv \qexp{(\Phip w_+ \Lambda,\eta) - (\Phip
\Lambda,\gamma)} a \cw{\Lambda} \mod{I^+_{w_+}} \\
\ctw{\Lambda} a \equiv \qexp{  (\Phim
\Lambda,\gamma) - (\Phim w_- \Lambda,\eta)} a \ctw{\Lambda}
\mod{I^-_{w_-}}.
\end{gather*}
\end{lem}

\begin{pf} The first identity  follows from Corollary \ref{cor 3.10} and the
definition of $I^+_{w_+}$. The second  identity can be deduced from
the first one by applying $\sigma$. 
\end{pf}

\smallskip

We continue to denote by $\cw{\Lambda} $ and $\ctw{\Lambda}$ the
images of these elements in $A / I_w$. It follows from Lemma \ref{lem
4.2} that the sets
$$\cal{E}_{w_+} = \{ \alpha \cw{\Lambda} \, \mid \, \alpha \in \C^*,
\Lambda \in \bL^+\}, \ \cal{E}_{w_-} = \{\alpha \ctw{\Lambda} \, \mid
\, \alpha \in \C^*, \Lambda \in \bL^+\}, \  \cal{E}_w = \cal{E}_{w_+}
\cal{E}_{w_-}$$
are multiplicatively closed sets of normal elements in $A/I_w$. Thus
$\cal{E}_w$ is an Ore set in $A/I_w$.
Define
$$A_w = \left( A/I_w \right)_{\cal{E}_w}.$$
Notice that $\sigma$ extends to an isomorphism $\sigma :
\hat{A}_{\hat{w}} \to A_w$, where $\hat{w} = (w_-,w_+).$ 

\smallskip

\begin{prop} \label{prop 4.3}
For all $w \in W \times W$, $A_w \neq (0)$. 
\end{prop}

\begin{pf} Notice first that since the generators of $A_w$ and the
elements of $\cal{E}_w$ are $\bL \times \bL$ homogeneous, it suffices
to work in the one-parameter case.   
The proof is then similar to  that of \cite[Theorem 2.2.3]{HL1}
(written in the $SL(n)$-case). For
completeness we recall the steps of this proof. The technical details
are straightforward generalizations to the general case of \cite[loc.
cit.]{HL1}.

 For  $1 \leq i \le n$ denote by $U_q(\frak{sl}_i(2))$ the Hopf
subalgebra of $\Uq$ generated by $e_i,f_i,k_i^{\pm 1}.$  The
associated quantized function algebra $A_i \cong \C_q[SL(2)]$ is 
naturally a quotient of $A$. Let $\sigma_i$ be the reflection
associated to the root $\alpha_i$. It is easily seen that there exist
$M_i^+$ and $M_i^-$, non-zero $(A_i)_{(\sigma_i,e)}$ and
$(A_i)_{(e,\sigma_i)}$ modules respectively. These modules can then be
viewed as non-zero $A$-modules.

  Let $w_+ = \sigma_{i_1} \dots \sigma_{i_k}$ and $w_- =
\sigma_{j_1} \dots \sigma_{j_m}$ be reduced expressions for $w_{\pm}$.
Then 
$$M^+_{i_1} \otimes \dots \otimes  M^+_{i_k} \otimes M^-_{j_1} \otimes
\dots \otimes M^-_{j_m} $$
is a non-zero $A_w$-module. 
\end{pf}

\smallskip

In the one-parameter case the proof of the following result was found
independently by the authors in \cite[1.2]{HL2} and Joseph in
\cite[6.2]{J1}.

\smallskip

\begin{thm} \label{thm 4.4} 
Let $P \in \Spec \Cqp$. There exists a unique $w \in W \times W$ such
that $P \supset I_w$ and $(P/I_w) \cap \cal{E}_w =\emptyset$.
\end{thm}

\begin{pf} Fix a dominant weight $\Lambda$. Define an ordering on the
weight vectors of $L(\Lambda)^*$ by $f \le f'$ if $f' \in \uqpbp f.$ 
This is a preordering which induces a partial ordering on the set of
one dimensional weight spaces. Consider the set:
$$\cal{F}(\Lambda) = \{ f \in L(\Lambda)^*_\mu \, \mid \,
\cfv{f}{v_\Lambda} \notin P \}.$$
Let $f$ be an element of $\cal{F}(\Lambda)$ which is maximal for the
above ordering. Suppose that $f'$ has the same property and that $f$
and $f'$ have weights $\mu$ and $\mu'$ respectively. By Corollary
\ref{cor 3.10} the two elements $\cfv{f}{v_\Lambda} $ and
$\cfv{f'}{v_\Lambda}$ are normal modulo $P$. Therefore we have, modulo
$P$,
$$\cfv{f}{v_\Lambda} \cfv{f'}{v_\Lambda} = \qexp{(\Phip \Lambda,\Lambda)
- (\Phip \mu,\mu')} \cfv{f'}{v_\Lambda} \cfv{f}{v_\Lambda} =
\qexp{2(\Phip \Lambda,\Lambda) - (\Phip \mu,\mu') - (\Phip \mu',\mu)}
\cfv{f}{v_\Lambda} \cfv{f'}{v_\Lambda}.$$
But, since $u$ is alternating, $2(\Phip \Lambda,\Lambda) - (\Phip
\mu,\mu') - (\Phip \mu',\mu) = 2(\Lambda,\Lambda) - 2(\mu,\mu').$
Since $P$ is prime and $q$ is not a root of unity we can deduce that
$(\Lambda,\Lambda) = (\mu,\mu').$ This forces $\mu = \mu' \in
W(-\Lambda)$. In conclusion, we have shown that for all dominant
$\Lambda$ there exists a unique (up to scalar multiplication) maximal
element $g_\Lambda \in \cal{F}(\Lambda)$ with weight $-w_\Lambda
\Lambda,$ $w_\Lambda \in W$. Applying the argument above to a pair of
such elements, $c_{g_\Lambda, v_\Lambda} $ and
$c_{g_\Lambda,v_{\Lambda'}}$, yields that $(w_\Lambda
\Lambda,w_{\Lambda'} \Lambda') = (\Lambda,\Lambda')$ for all $\Lambda,
\Lambda' \in \bL^+$. Then it is not difficult to show that this
furnishes a unique $w_+ \in W$ such that $w_+ \Lambda = w_\Lambda
\Lambda$ for all $\Lambda \in \bL^+$. Thus for each $\Lambda \in
\bL^+$,
$$c_{g,v_\Lambda} \in P \Longleftrightarrow g \nleq f_{-w_+\Lambda}.$$
Hence $P \supset I^+_{w_+}$ and $P \cap \cal{E}_{w_+} = \emptyset.$
It is easily checked that such a $w_+$ must be unique. Using $\sigma$
one deduces the existence and uniqueness of $w_-$.
\end{pf}

\smallskip

\begin{defn} A prime ideal $P$ such that $P \supset I_w$ and $P
\cap \cal{E}_w = \emptyset$ will be called a prime ideal of type
$w$. We denote by $\Spec_w \Cqp$, resp. $\Prim_w \Cqp$, the subset of
$\Spec \Cqp$ consisting of prime, resp. primitive, ideals of type $w$.
\end{defn}

\smallskip

Clearly $\Spec_w \Cqp \cong \Spec A_w$  and $\sigma(\Spec_{\hat{w}}
\C_{q^{-1},p}[G]) = \Spec_w \Cqp$. The following corollary is
therefore clear.

\smallskip

\begin{cor} \label{cor 4.5}
One has
 $$\Spec \Cqp = \sqcup_{w \in W \times W} \Spec_w \Cqp, \quad \Prim \Cqp =
\sqcup_{w \in W \times W} \Prim_w \Cqp.$$ 
\end{cor}

\smallskip

We end this section by a result which is the key idea in \cite{J1} for
analyzing the adjoint action of $A$ on $A_w$.  It says that in the one
parameter case the quantized function algebra $\C_q[B^-]$ identifies
with $\uqbp$ through the  Rosso-Tanisaki-Killing form,
\cite{DCL,HL3,J1}. Evidently this continues to hold in the
multi-parameter case. For completeness we include a proof of that
result.

Set $\C_{q,p}[B^-] = A/I_{(w_0,e)}.$ The
embedding $\uqpbm \to \Dqp$ induces a Hopf algebra map $\phi : A \to
\uqpbm^{\circ}$, where $\uqpbm^{\circ}$ denotes the cofinite dual. On
the other hand the non-degenerate Hopf algebra pairing $\hpair{\,}{\,}_{p^{-1}}$
furnishes an injective morphism $\theta: \uqpbp^{\text{op}} \to \uqpbm^*$.

\smallskip

\begin{prop} \label{prop 4.6}
\rom{1.} $\C_{q,p}[B^-]$ is an $\bL$-bigraded Hopf algebra.  

\rom{2.} The map $\gamma = \theta^{-1} \phi: \C_{q,p}[B^-] \to
\uqpbp^{\text{op}}$ is an   isomorphism of Hopf algebras. 
\end{prop}

\begin{pf}
1. It is easy to check that $I_{(w_0,e)}$ is an $\bL \times \bL$
graded bi-ideal of the bialgebra $A$. Let $\mu \in \bL^+$ and fix a dual
basis $\{v_\nu; f_{-\nu}\}_\nu$ of $L(\mu)$ (with the usual abuse of notation).  Then 
$$\sum_\nu \cfv{v_\nu}{f_{-\eta}}
\cfv{f_{-\nu}}{v_\gamma} = \sum_\nu S(\cfv{f_{-\eta}}{v_\nu})
\cfv{f_{-\nu}}{v_\gamma} =
\epsilon(\cfv{f_{-\eta}}{v_\gamma}).$$
Taking $\gamma = \eta = \mu$ yields $\tilde{c}_\mu c_\mu = 1$ modulo
$I_{(w_0,e)}$. If $\gamma = w_0\mu$ and $\eta \neq w_0\mu$, the above
relation shows that $S(\cfv{f_{-\eta}}{v_{w_0\mu}})
\tilde{c}_{-w_0\mu} \in I_{(w_0,e)}$. Thus $I_{(w_0,e)}$ is a Hopf ideal.

2. We first show that 
\begin{equation} \label{surj}
\forall \, \Lambda \in \bL^+, \cfv{f}{v_\Lambda} \in
C(\Lambda)_{-\lambda,\Lambda}, \ \exists ! \, x_\lambda \in
U^+_{\Lambda - \lambda}, \quad \phi(\cfv{f}{v_\Lambda}) =
\theta(x_\lambda \cdot k_{-\Lambda}).
\end{equation}
Set $c=\cfv{f}{v_\Lambda}$. Then $c(U^-_{-\eta}) = 0$ unless $\eta =
\Lambda -\lambda$; denote by $\bar{c}$ the restriction of $c$ to
$U^-$. By the non-degeneracy of the pairing on
$U^+_{\Lambda -\lambda} \times U^-_{\lambda - \Lambda}$ we know that
there exists a unique $x_\lambda \in U^+_{\Lambda - \lambda}$ such that
$\bar{c} = \theta(x_\lambda)$. Then, for all $y \in
U^-_{\lambda - \Lambda}, $ we have
\begin{align*}
c(y \cdot k_\mu) &= \pair{f}{y \cdot k_\mu \cdot v_\Lambda} =
\qexp{-(\Phim \Lambda,\mu)} \bar{c}(y) 
= \qexp{-(\Phim \Lambda,\mu)}\hpair{x_\lambda}{y} \\
& = \hpair{x_\lambda \cdot
k_{-\Lambda}}{y \cdot k_\mu}_{p^{-1}}
\end{align*}
by \eqref{pairing}. This proves \eqref{surj}.

 We now show that $\phi$ is injective on $A^+$. Suppose that
$c = \cfv{f}{v_\Lambda} \in C(\Lambda)_{-\lambda,\Lambda} \cap \Ker
\phi$, hence $c =  0$ on $\uqpbm$. Since $L(\Lambda) = \uqpbm v_\Lambda
= \Dqp v_\Lambda$ it follows that $c= 0$.  An easy weight argument using
\eqref{surj} shows then that  $\phi$ is injective on
$A^+$.


It is clear that $\Ker \phi \supset I_{(w_0,e)}$, and that $A^+ A^- =
A$ implies $\phi(A) = \phi(A^+[\tilde{c}_\mu \, ; \, \mu \in \bL^+]).$
Since $\tilde{c}_\mu = c_\mu^{-1}$ modulo $I_{(w_0,e)}$ by part 1, if
$a \in A$ there exists $\nu \in \bL^+$ such that $\phi(c_\nu) \phi(a)
\in \phi(A^+)$. The inclusion  $\Ker \phi \subset I_{(w_0,e)}$ follows
easily. Therefore $\gamma$ is a well defined Hopf algebra
morphism. 

If  $\alpha_j \in \bB$,
there exists $\Lambda \in \bL^+$ such that $L(\Lambda)_{\Lambda -
\alpha_j} \ne (0)$. Pick $0 \ne f \in L(\Lambda)^*_{-\Lambda +
\alpha_j}$. Then \eqref{surj} shows that, up to some scalar,
$\phi(\cfv{f}{v_\Lambda}) = \theta( e_j \cdot k_{-\Lambda})$. 
If  $\lambda \in \bL$, there  exists $ \Lambda \in W\lambda \cap
\bL^+$; in particular $L(\Lambda)_{\lambda} \ne (0).$ Let $v \in
L(\Lambda)_{\lambda}$ and $ f \in L(\Lambda)^*_{-\lambda}$ such
that $\pair{f}{v} =1$. Then it is easily verified that
$\phi(\cfv{f}{v}) = \theta(k_{-\lambda}).$ This proves that $\gamma$
is surjective, and the proposition.
\end{pf}

\smallskip

\subsection{The adjoint action}
Recall that if $M$ is an arbitrary $A$-bimodule one defines the
adjoint action of $A$ on $M$ by
$$\forall \, a \in A, \, x \in M, \quad \ad(a).x = \sum a_{(1)} x
S(a_{(2)}).$$ 
Then it is well known that the subspace of $\ad$-invariant elements
$M^{\text{ad}} = \{x \in M \, \mid \,  \forall a \in A, \, \ad(a).x =
\epsilon(a)x\}$ is equal to $\{x \in M \, \mid \, \forall a \in A, \,
ax = xa \}.$ 

Henceforth we fix $w \in W \times W$ and work inside $A_w$. For
$\Lambda \in \bL^+$, $f \in L(\Lambda)^*$ and $v \in L(\Lambda)$ we
set
$$z^+_f = \cw{\Lambda}^{-1} \cfv{f}{v_\Lambda}, \quad z^-_v =
\ctw{\Lambda}^{-1} \cfv{v}{f_{-\Lambda}}.$$ 
 Let $\{\omega_1, \dots,\omega_n\}$ be a basis of $\bL$ such that
$\omega_i \in \bL^+$ for all $i$. Observe that $\cw{\Lambda}
\cw{\Lambda'} $ and $\cw{\Lambda'}\cw{\Lambda}$ differ by a non-zero
scalar (similarly for $\ctw{\Lambda}\ctw{\Lambda'}  $). For each
$\lambda = \sum_i \ell_i \omega_i \in \bL$ we define normal 
elements of $A_w$ by 
$$\cw{\lambda} = \prod_{i=1}^n \cw{ \omega_i}^{\ell_i}, \quad
\ctw{\lambda} = \prod_{i=1}^n \ctw{ \omega_i}^{\ell_i},
 \quad
d_\lambda = (\ctw{\lambda} \cw{\lambda})^{-1}.$$
Notice then that, for $\Lambda \in \bL^+$,  the ``new''
$\cw{\Lambda}$ belongs to 
$\C^* \cfv{f_{-w_+\Lambda}}{v_\Lambda}$ (similarly for
$\ctw{\Lambda}$). 
 Define subalgebras of $A_w$ by
\begin{gather*}
C_w = \C[z^+_f, z^-_v, \cw{\lambda} \, ; \, f \in L(\Lambda)^*, v \in
L(\Lambda),    \, \Lambda \in \bL^+, \, \lambda \in \bL] \\
C^+_w = \C[z^+_f \, ; f \in L(\Lambda)^*, \,\Lambda \in \bL^+],
\quad  C^-_w = \C[z^-_v \, ; v \in L(\Lambda), \,\Lambda \in \bL^+].
\end{gather*}
Recall that the torus $H$ acts on $A_{\lambda,\mu}$ by $r_h(a) =
\mu(h)a$, where $\mu(h) = \pair{\mu}{h}$. 
Since the generators of $I_w$ and the elements of
$\cal{E}_w$ are eigenvectors for $H$, the action of $H$ extends to an
action on $A_w$. The algebras $C_w$ and $C^\pm_w$ are obviously
$H$-stable.

\smallskip

\begin{thm}\label{thm 4.7}
\rom{1}. $C_w^H = \C[z^+_f, z^-_v \, ; \, f \in L(\Lambda)^*, v \in
L(\Lambda), \, \Lambda \in \bL^+].$

\rom{2.} The set $\cal{D} = \{d_\Lambda \, ; \, \Lambda \in \bL^+\}$
is an Ore subset of $C_w^H$. Furthermore $A_w = (C_w)_{\cal{D}}$ and
$A_w^H = (C_w^H)_{\cal{D}}.$

\rom{3.} For each $\lambda \in \bL$,  
let $(A_w)_\lambda = \{ a \in A_w \, \mid  \, r_h(a) = \lambda(h)a \}$. Then
 $A_w = \bigoplus_{\lambda \in \bL} (A_w)_\lambda$ and 
$(A_w)_\lambda = A_w^H c_{w\lambda}$ . Moreover each $(A_w)_\lambda$
is an ad-invariant subspace. 
\end{thm}

\begin{pf} 
Assertion 1 follows from
$$\forall h \in H, \quad r_h(z^\pm_f) = z^\pm_f, \quad
r_h(\cw{\lambda}) = \lambda(h) \cw{\lambda}, \quad r_h(\ctw{\lambda}) =
\lambda(h)^{-1} \ctw{\lambda}.$$

 Let $\{v_i ; f_i\}_i$ be a dual basis for $L(\Lambda)$. Then
$$1 =\epsilon(\cfv{f_{-\Lambda}}{v_\Lambda})      = \sum_i
S(\cfv{f_{-\Lambda}}{v_i})\cfv{f_i}{v_\Lambda}  = \sum_i
\cfv{v_i}{f_{-\Lambda}}\cfv{f_i}{v_\Lambda}.$$
Multiplying both sides
of the equation by $d_\Lambda$ and using the normality of $c_{w\Lambda}$
and $\ct_{w\Lambda}$ yields $d_\Lambda = \sum_i a_i z^-_{v_i}
z^+_{f_{i}}$ for some $a_i \in \C$. Thus $\cal{D} \subset C^H_w$. Now
by Theorem \ref{thm 4.1}  any
element of $A_w$ can be written in the form
$\cfv{f_1}{v_1}\cfv{f_2}{v_2}d_\Lambda^{-1}$ where $v_1 = v_{\Lambda_1}$, $v_2
= v_{-\Lambda_2}$ and $\Lambda_1, \Lambda_2, \Lambda \in
\bL^+$.
 This element lies in $(A_w)_\lambda$ if and only if $\Lambda_1
-\Lambda_2 = \lambda$. In this case
$\cfv{f_1}{v_1}\cfv{f_2}{v_2}d_\Lambda^{-1}$ is equal, up to a scalar,
to the element $z^+_{f_1}z^-_{f_2}d_{\Lambda+\Lambda_2}^{-1}
\cw{\lambda} \in (C^H_w)_{\cal{D}} \cw{\lambda}$. Since the adjoint
action commutes with the right action of $H$, $(A_w)_\lambda$ is an
ad-invariant subspace. The remaining assertions then follow.
\end{pf}

\smallskip

	We now study the adjoint action of $\cqp$ on $A_w$. The key result is
Theorem \ref{thm 4.12}.

\smallskip

\begin{lem}\label{lem 4.8}
Let $T_\Lambda = \{ z^+_f \mid f \in L(\Lambda)^* \}$. Then 
 $C^+_w = \bigcup_{\Lambda \in \bL} T_\Lambda$.
\end{lem}

\begin{pf}
It suffices to prove that if $\Lambda, \Lambda' \in \bL^+$ and $f \in
L(\Lambda)^*$, then there exists a $g \in L(\Lambda + \Lambda')^*$
such that $z^+_f = z^+_g$. Clearly we may assume that $f$ is a weight
vector. Let $\iota: L(\Lambda+\Lambda') \to L(\Lambda)\otimes
L(\Lambda')$ be the canonical map. Then
$$c_{f,v_\Lambda}c_{f_{-w_+\Lambda'}, v_{\Lambda'}} =
c_{f_{-w_+\Lambda'} \otimes f, v_\Lambda \otimes v_{\Lambda'}} =
c_{g,v_{\Lambda +\Lambda'}}$$ where $g = \iota^*(f_{-w_+\Lambda'}
\otimes f)$. Multiplying the images of these elements in $A_w$ by the
inverse of $c_{w(\Lambda+\Lambda')} \in \C^* c_{w\Lambda}c_{w\Lambda'}$
yields the desired result.
\end{pf}

\smallskip

\begin{prop} \label{prop 4.9}
 Let $E$ be an object of $\cal{C}_{q,p}$ and let $\Lambda
\in \bL^+$.  Let $\sigma:
L(\Lambda) \to E \otimes L(\Lambda) \otimes E^*$ be  
the map $(1\otimes \psi^{-1})(\iota\otimes 1)$ where $\iota : \C \to E
\otimes E^*$ is  the canonical
embedding and $\psi^{-1}: E^* \otimes L(\Lambda) \to L(\Lambda) \otimes E^*$ is the
inverse of the braiding map described in \S 3.5.
 Then for any
$ c =c_{g,v} \in C(E)_{-\eta,\gamma} $ and $f \in L(\Lambda)^*$ 
$$ 	\ad(c) .z^+_f = q^{(\Phi_+ w_+ \Lambda, \eta)} \,  z^+_{\sigma^*(v\otimes f
\otimes g)}.  $$
In particular $C_w^+$ is a locally finite $\Cqp$-module for the adjoint action.
\end{prop}

\begin{pf} Let $\{v_i;g_i\}_i$ be a dual basis of $E$ where
$v_i \in E_{\nu_i}, \, g_i \in E^*_{-\nu_i} $. Then $\iota(1) = \sum v_i \otimes 
g_i$. By \eqref{braiding} we have
$$\psi^{-1} (g_i \otimes v_\Lambda) = a_i (v_\Lambda \otimes g_i)$$
where $a_i = q^{-(\Phi_+\Lambda,\nu_i)} = q^{(\Phi_- \nu_i, \Lambda)}$.  
On the other hand the commutation relations given in Corollary \ref{cor 3.10} imply that $c_{g,v_i} 
c_{w\Lambda}^{-1} =
b a_i c_{w\Lambda}^{-1} c_{g,v_i}$,
where $b = q^{(\Phi_+w_+ \Lambda,\eta)}$.  Therefore
$$ 
\ad(c) .z^+_f = \sum b
a_i c_{w\Lambda}^{-1} c_{g,v_i} c_{f,v_\Lambda} c_{v,g_i} = b
c_{w\Lambda}^{-1} c_{v \otimes f \otimes g, \sum a_i v_i \otimes
v_\Lambda \otimes g_i} = b c_{w\Lambda}^{-1} c_{v \otimes f \otimes g,
\sigma(v_\Lambda)} .
$$ 
Since the map $\sigma$ is a morphism of $\Dqp$-modules it is easy to
see that $c_{v \otimes f \otimes g, \sigma(v_\Lambda)} = c_{\sigma^*(v
\otimes f \otimes g), v_\Lambda}.$
\end{pf}

\smallskip

Let $\gamma: \Cqp \to \uqpbp$ be the algebra anti-isomorphism
 given in  Proposition \ref{prop 4.6}.

\smallskip 

\begin{lem} \label{lem 4.10}
 Let $ c  = c_{g,v} \in  \cqp_{-\eta,\gamma}$, $f \in L(\Lambda)^*$ be
as in the previous theorem and $ x \in \uqpbp$  be such that
$\gamma(c) = x$. Then 
$$c_{S^{-1}(x).f,v_\Lambda}
=c_{\sigma^*(v\otimes f \otimes g), v_\Lambda}.$$
\end{lem}

\begin{pf} Notice that it suffices to show that 
$$c_{S^{-1}(x).f,v_\Lambda}(y) =c_{\sigma^*(v\otimes f \otimes g), 
v_\Lambda}(y)$$
for all $y \in \uqpbm$. Denote by $\langle \ \mid \ \rangle$ the Hopf pairing
$\langle \ \mid \ \rangle_{p^{-1}}$ between $\uqpbp^{\text{op}}$ and
$\uqpbm$ as in \S \ref{3.4}.
Let $\chi$ be the one dimensional representation of $\uqpbp$ associated to
$v_\Lambda$ and let $ \tilde{\chi} = \chi \cdot \gamma$. Notice that 
$\chi(x) = \langle x \mid
t_{-\Lambda}\rangle$; so $\tilde{\chi}(c) = c(t_{-\Lambda})$.
Recalling that $\gamma$ is a morphism of coalgebras and using the
relation (c$_{xy})$ of \S \ref{2.3} in the double $\uqpbp \Join
\uqpbm$, we obtain
\begin{eqnarray*}
c_{S^{-1}(x).f, v_\Lambda}(y) & = & f(xyv_\Lambda) \\
& = &\sum \langle x_{(1)}\mid y_{(1)}\rangle \langle x_{(3)}\mid
S(y_{(3)})\rangle  \,
f(y_{(2)}x_{(2)}v_\Lambda)  \\
&= &\sum \langle x_{(1)}\mid y_{(1)} \rangle \langle x_{(3)}\mid
S(y_{(3)})\rangle  \,
\chi(x_{(2)}) \,
f(y_{(2)}v_\Lambda)  \\
&=& \sum \langle x_{(1)}\chi(x_{(2)})\mid y_{(1)}\rangle \langle
x_{(3)}\mid S(y_{(3)}\rangle \, 
f(y_{(2)}v_\Lambda)  \\
&=&\sum  (c_{(1)}\tilde{\chi}(c_{(2)}))(y_{(1)}) \,  c_{(3)} 
(S(y_{(3)})) \,
f(y_{(2)}v_\Lambda)\\
&=&\sum r_{\tilde{\chi}}(c_{(1)})(y_{(1)}) \,
c_{f,v_\Lambda}(y_{(2)}) \,
S(c_{(2)})(y_{(3)}). 
\end{eqnarray*}
Since
$r_{\tilde{\chi}}(c_{g,v_i}) = q^{(\Phi_-\nu_i,\Lambda)} c_{g,v_i}$,
one shows as in the proof of Proposition \ref{prop 4.9}  that
\begin{eqnarray*}
c_{S^{-1}(x).f, v_\Lambda}(y) &=& \sum
r_{\tilde{\chi}}(c_{(1)})(y_{(1)}) \, 
c_{f,v_\Lambda}(y_{(2)}) \, S(c_{(2)})(y_{(3)}) \\
&=& \sum q^{(\Phi_-\nu_i,\Lambda)} \, (c_{g,v_i} c_{f,v_\Lambda}
c_{v,g_i})(y) \\
&=& c_{\sigma^*(v\otimes f \otimes g),v_\Lambda}(y),
\end{eqnarray*}
as required.
\end{pf}

\smallskip

\begin{thm} \label{thm 4.11}
Consider $C^+_w$ as a  $\cqp$-module via the adjoint action. Then
\begin{enumerate}
\item  $\Soc C^+_w = \C$.
\item $\Ann C^+_w \supset I_{(w_0,e)}$.
\item The elements $c_{f_{-\mu}, v_\mu}$, $\mu \in \bL^+$, 
 act diagonalizably on $C^+_w$.
\item $\Soc C^+_w = \{z \in C^+_w \mid \Ann z \supset I_{(e,e)}\}$.
\end{enumerate}
\end{thm}

\begin{pf}
For  $\Lambda \in \bL^+,$ define a  $\uqpbp$-module by
$$S_\Lambda = (\uqpbp
v_{w_+\Lambda})^* = L(\Lambda)^*/(\uqpbp v_{w_+\Lambda})^\perp.$$ 
It is easily checked that $\Soc S_\Lambda = \C f_{-w_+\Lambda}$ (see
\cite[7.3]{J1}).   Let $\delta: S_\Lambda \to T_\Lambda$ be
the linear map given by $\bar{f} \mapsto z^+_f$. Denote by $\zeta$ the
one-dimensional representation of $\cqp$ given by $\zeta(c) =
c(t_{-w_+\Lambda})$.  Let $c=c_{g,v} \in C(E)_{-\eta,\gamma}$.  Then
$l_\zeta(c) = q^{(\Phi_-\eta, w_+\Lambda)} c =
q^{-(\Phi_+w_+\Lambda,\eta))} c$.  Then, using Proposition \ref{prop
4.9} and Lemma \ref{lem 4.10} we obtain,
$$\ad(l_\zeta(c)).\delta(\bar{f}) = q^{-(\Phi_+w_+\Lambda, \eta)}\ad(c).
z^+_f = z^+_{S^{-1}\gamma(c).f} = \delta(S^{-1}(\gamma(c))\bar{f}).$$
Hence, $\ad(l_\zeta(c)).\delta(\bar{f}) =
\delta(S^{-1}(\gamma(c))\bar{f})$ for all $c \in A$. This immediately
implies part (2) since $\Ker \gamma
\supset I_{(w_0,e)}$ and $l_\zeta(I_{(w_0,e)}) = I_{(w_0,e)}$. If
$S_\Lambda$ is given the structure of an $A$-module via $S^{- 
1}\gamma$, then $\delta$ is a homomorphism from $S_\Lambda$ to the
module $T_\Lambda$ twisted by the automorphism $l_\zeta$.  Since
$\delta(f_{-w_+\Lambda}) = 1$ it follows that $\delta$ is bijective and
that $\Soc T_\Lambda = \delta(\Soc S_\Lambda) = \C$. Part (1) then
follows from Lemma \ref{lem 4.8}.  	Part (3) follows from  the
above formula and the fact that $\gamma(c_{f_{-\mu}, v_\mu}) =
s_{-\mu}$. Since $A/I_{(e,e)}$ is generated by the images of
the elements $c_{f_{-\mu}, v_\mu}$, 
(4) is a consequence of the definitions. 
\end{pf}

\smallskip

\begin{thm}\label{thm 4.12} 
Consider $C^H_w$ as a $\cqp$-module via the adjoint 
action. Then 
$$\Soc C^H_w = \C.$$
\end{thm}

\begin{pf}
By Theorem \ref{thm 4.11} we have that $\Soc C^+_w = \C$. Using the
map $\sigma$, one obtains analogous results for $C^-_w$. The map
$C^+_w \otimes C^-_w \to C^H_w$ is a module map for the adjoint
action which is surjective by Theorem
\ref{thm 4.1}. So it suffices 
to show that
$\Soc C^+_w \otimes C^-_w =\C$. The following argument is taken from 
\cite{J1}. 

By
the analog of Theorem \ref{thm 4.11} for $C^-_w$ we have that the
elements $c_{f_{-
\Lambda},v_\Lambda}$ act as commuting diagonalizable operators on $C^-_w$.
Therefore an element of $C^+_w \otimes C^-_w$ may be written as $\sum a_i 
\otimes
b_i$ where the $b_i$ are linearly independent weight vectors.
Let $\cfv{f}{v_\Lambda}$ be a generator of $I^+_e$. 
Suppose that $\sum a_i \otimes b_i \in \Soc (C^+_w \otimes C^-_w)$.
Then
\begin{equation*}
\begin{split}
0 = \ad(\cfv{f}{v_\Lambda}). (\sum_i a_i \otimes b_i ) &=
\sum_{i,j}\ad(\cfv{f}{v_j}) .a_i \otimes
\ad(\cfv{f_j}{v_\Lambda}).b_i \\
&= \sum_i \ad(\cfv{f}{v_\Lambda}).a_i \otimes
\ad(\cfv{f_{-\Lambda}}{v_\Lambda}).b_i \\
&=\sum_i \ad(\cfv{f}{v_\Lambda}).a_i\otimes \alpha_i b_i
\end{split}
\end{equation*}
for some $\alpha_i \in \C^*$.  Thus
$\ad(\cfv{f}{v_\Lambda}).a_i = 0$ for all $i$. Thus the $a_i$ are annihilated by 
the left ideal
generated by the  $\cfv{f}{v_\Lambda}$. But this left ideal is two-sided modulo
$I_{(w_0,e)}$ and $\Ann C^+_w \supset I_{(w_0,e)}$. 
 Thus the $a_i$ are annihilated by $I_{(e,e)}$ and so lie in $\Soc C^+_w$ by 
Theorem \ref{thm 4.11}. 
Thus $\sum a_i \otimes b_i \in  \Soc  (\C \otimes C^-_w) = \C \otimes
\C.$ \end{pf}

\smallskip

\begin{cor} \label{cor 4.13}
 The algebra $A_w^H$ contains no nontrivial ad-invariant ideals.
Furthermore,  $(A^H_w)^{\ad} = \C$.
\end{cor}

\begin{pf} Notice that Theorem \ref{thm 4.12} implies that $C^H_w$
contains no nontrivial ad-invariant ideals. 
Since  $A_w^H$ is a localization of $C_w^H$ the same must 
be true for $A_w^H$.
 Let $a \in (A^H_w)^{\ad}\backslash \C$. Then $a$ is central and so for any 
$\alpha \in
\C$,  $(a-\alpha)$ is a non-zero ad-invariant ideal of $A^H_w$.  This 
implies that $a-
\alpha$ is invertible in $A^H_w$ for any $\alpha \in \C$. This contradicts 
the fact that
$A^H_w$ has countable dimension over $\C$.
\end{pf}

\smallskip

\begin{thm} \label{thm 4.14}
Let $Z_w$ be the center of $A_w$. Then 
\begin{enumerate}
\item $Z_w = A_w^{ad}$;
\item $Z_w = \bigoplus_{\lambda\in \bL}Z_\lambda$ where $Z_\lambda = Z_w \cap
A^H_w c_{w\lambda}$;
\item If $Z_\lambda \neq (0)$, then $Z_\lambda = \C u_\lambda$ for some unit
$u_\lambda$; 
\item The group $H$ acts transitively on the maximal ideals of $Z_w$.
\end{enumerate}
\end{thm}

\begin{pf}The proof of (1) is standard. Assertion (2) follows from
Theorem \ref{thm 4.7}.
 Let $u_\lambda$ be a non-zero element of $Z_\lambda$. Then $u_\lambda = 
 a c_{w\lambda}$, for some $a\in A^H_w$. 
 This implies that $a$ is normal and hence $a$ generates an ad-invariant 
 ideal of $A^H_w$. Thus $a$ (and hence also $u_\lambda$) is a unit by 
Theorem \ref{cor 4.13}. Since
$Z_0 = \C$, it follows that $Z_\lambda = \C u_\lambda$. Since 
 the action of $H$ is given by $r_h(u_\lambda) = \lambda(h) u_\lambda$, it is 
clear that $H$
acts transitively on the maximal ideals of $Z_w$.
\end{pf}

\smallskip

\begin{thm} \label{thm 4.15}
The ideals of $A_w$ are generated by their intersection with the center, $Z_w$.
\end{thm}

\begin{pf}
Any element $f \in A_w$ may be written uniquely in the form $f = \sum a_\lambda
c_{w\lambda}$ where $a_\lambda \in A_w^H$. Define  $\pi : A_w \rightarrow 
A_w^H$ to be
the projection given by $\pi(\sum a_\lambda c_{w\lambda}) = a_0$ and notice 
that $\pi$ is a
module map for the adjoint action.  Define the support of $f$ to be 
$Supp(f) = \{ \lambda \in \bL \mid a_\lambda \neq 0\}$.
Let $I$ be an ideal of $A_w$. For any set $Y \subseteq \bL$ such that $0 \in 
Y$ define
$$I_Y = \{ b \in A_w^H \mid b = \pi(f) \mbox{ for some } f \in I \mbox{ such 
that }
Supp(f) \subseteq Y\}$$
 If $I$ is ad-invariant then $I_Y$ is an ad-invariant ideal of $A_w^H$ and 
hence is either
$(0)$ or $A_w^H$.

Now let $I' = (I \cap Z_w)A_w$ and suppose that $I\neq I'$. Choose an 
element $f=\sum
a_\lambda c_{w\lambda} \in I \backslash I'$ whose support $S$ has the smallest 
cardinality.
We may assume without loss of generality that $0 \in S$. Suppose that there 
exists $g \in
I'$ with $Supp(g) \subset S$. Then there exists a $g'\in I'$ with $Supp(g') 
\subset S$ and
$\pi(g') = 1$. But then $f - a_0g'$ is an element of $I'$ with smaller 
support than $F$.
Thus there can be no elements in $I'$ whose support is contained in $S$. So 
we may
assume that $\pi(f)=a_0 =1$. For any $c \in \cqp$, set $f_c = \ad(c).f - 
\epsilon(c)f$.
Since $\pi(f_c) =0$ it follows that 
$|Supp(f_c)| < |Supp(f)|$ and hence that $f_c=0$. Thus $f \in I\cap 
A_w^{\ad}= I\cap
Z_w$, a contradiction.
\end{pf}

\smallskip

	Putting these results together yields the  main theorem of
this section, which completes Corollary \ref{cor 4.5} by describing
the set of primitive ideals of type $w$.

\smallskip

\begin{thm} \label{thm 4.16}
For  $w \in W \times W$ the subsets $\Prim_w \cqp$ are precisely the
$H$-orbits inside $\Prim \cqp$.
\end{thm}

\smallskip

	Finally we calculate the  size of these orbits in the algebraic 
case. 
Set $\bL_w = \{ \lambda \in \bL \, \mid \, Z_\lambda \neq (0)\}$.
Recall the definition of $s(w)$ from \eqref{s(w)} and that $p$ is
called $q$-rational if $u$ is algebraic. In this case we know by
Theorem \ref{thm 1.8} that there exists $m \in \N^*$ such that
$\Phi(m\bL) \subset \bL$.

\smallskip
 
\begin{prop} \label{prop 4.17} 
Suppose that $p$ is $q$-rational.  Let $\lambda \in \bL$ and
$y_\lambda = c_{w\Phi_- m\lambda}\tilde{c}_{w\Phi_+m\lambda}$. Then
\begin{enumerate}
\item $y_\lambda$ is ad-semi-invariant. In fact, for any $c \in 
A_{-\eta,\gamma}$,
$$      \ad(c) . y_\lambda = \qexp{(m \sigma(w)\lambda, 
\eta)}\epsilon(c)y_\lambda.$$
where $\sigma(w) = \Phi_-w_-\Phi_+ -\Phi_+w_+\Phi_-$
\item $\bL_w \cap 2m\bL = 2\Ker \sigma(w) \cap m\bL$
\item $\dim Z_w = n - s(w)$
\end{enumerate}
\end{prop}

\begin{pf}
Using Lemma \ref{lem 4.2}, we have that for $c \in 
A_{-
\eta,\gamma}$
\begin{align*}
cy_\lambda &= \qexp{(\Phi_+w_+\Phi_-m\lambda,-\eta)}\qexp{(\Phi_+\Phi_-
m\lambda,\gamma)}
\qexp{(\Phi_-w_-\Phi_+m\lambda,\eta)}\qexp{(\Phi_-\Phi_+m\lambda,-
\gamma)}y_\lambda c \\
	&= \qexp{(m \sigma(w)\lambda, \eta)}y_\lambda c.
\end{align*}
From this it follows easily that
$$      \ad(c) . y_\lambda = \qexp{(m \sigma(w)\lambda, 
\eta)}\epsilon(c)y_\lambda.$$

Since (up to some scalar) 
$y_\lambda =
d^{-1}_{\Phi m \lambda} d^{-1}_{m \lambda}c^{-2}_{wm 
\lambda}$ it follows from
Theorem \ref{thm 4.7} that $ y_\lambda \in (A_w)_{-2m\lambda}$.
However, as a $\cqp$-module via the adjoint action, $A^H_w y_\lambda
\cong A^H_w
\otimes \C y_\lambda$ and hence $\Soc A^H_w y_\lambda = \C y_\lambda$.
Thus $Z_{-2m\lambda}
\neq (0)$ if and only if $y_\lambda$ is ad-invariant; that is, if and only if
$m \sigma(w) \lambda =0$. Hence
\begin{align*}
\dim Z_w &= \rk \bL_w = \rk (\bL_w \cap 2m\bL) = \rk \Ker_{m\bL} \sigma(w) \\
	&= \dim \Ker_{\h^*}\sigma(w) = n - s(w)
\end{align*}
as required.
\end{pf}

\smallskip

Finally, we may deduce that in the algebraic case the size of the of
the $H$-orbits $\Symp_w G$ and $\Prim_w \cqp$ are the same, cf.
Theorem \ref{thm 1.9}.

\smallskip

\begin{thm} \label{thm 4.18}
Suppose that $p$ is $q$-rational and let $w \in W \times W$. Then 
$$\forall P \in \Prim_w \cqp, \quad \dim (H/\Stab_H P)  = n - s(w).$$
\end{thm}
\begin{pf} This follows easily from  theorems \ref{thm 4.15},
\ref{thm 4.16} and Proposition \ref{prop 4.17}.
\end{pf}

\end{document}